\documentclass[a4paper,leqno,draft]{article}
\usepackage{amssymb,amsmath,amsfonts,amsthm,%eucal,
mathrsfs}

\newtheorem{theorem}{Theorem}[section]

\newtheorem{proposition}[theorem]{Proposition}

\theoremstyle{definition}
\newtheorem{remark}[theorem]{Remark}
\newtheorem{example}[theorem]{Example}

\newcommand{\wt}[1]{\widetilde{#1}}

\newcommand{\Cinf}{\ensuremath{\mathcal{C}^\infty}}
\newcommand{\Cinfc}{\ensuremath{\mathcal{C}^\infty_{\text{c}}}}
\newcommand{\D}{\ensuremath{{\cal D}}}
\renewcommand{\S}{\mathscr{S}}
\newcommand{\E}{\ensuremath{{\cal E}}}
\newcommand{\OM}{\ensuremath{{\cal O}_\mathrm{M}}}

\renewcommand{\L}{\mathcal{L}}

\newcommand{\mb}[1]{\ensuremath{\mathbb{#1}}}
\newcommand{\N}{\mb{N}}

\newcommand{\R}{\mb{R}}
\newcommand{\C}{\mb{C}}

\newcommand{\G}{\ensuremath{{\cal G}}}
\newcommand{\Gt}{\ensuremath{{\cal G}_\tau}}
\newcommand{\GtS}{{\mathcal{G}}_{\tau,\S}}
\newcommand{\Gc}{\ensuremath{{\cal G}_\mathrm{c}}}
\newcommand{\Gcinf}{\ensuremath{{\cal G}^\infty_\mathrm{c}}}
\newcommand{\GS}{\G_{{\, }\atop{\hskip-4pt\scriptstyle\S}}\!}
\newcommand{\EM}{\ensuremath{{\cal E}_{M}}}

\newcommand{\Et}{\ensuremath{{\cal E}_{\tau}}}
\newcommand{\ES}{\mathcal{E}_{\S}}
\newcommand{\EMinf}{\ensuremath{{\cal E}^\infty_{M}}}

\newcommand{\EMp}{\mathcal{E}_{M,p}}
\newcommand{\Nt}{\ensuremath{{\cal N}_{\tau}}}
\newcommand{\Neg}{\mathcal{N}}
\newcommand{\NS}{\mathcal{N}_{\S}}
\newcommand{\Npq}{\mathcal{N}_{p,q}}
\newcommand{\Npp}{\mathcal{N}_{p,p}}

\newcommand{\Gpq}{\mathcal{G}_{p,q}}
\newcommand{\Gpp}{\mathcal{G}_{p,p}}

\newcommand{\g}{\mathrm{g}} % for index `generalized'
\newcommand{\Ginf}{\ensuremath{\G^\infty}}

\newcommand{\GSinf}{\G^\infty_{{\, }\atop{\hskip-3pt\scriptstyle\S}}}

\newcommand{\lara}[1]{\langle #1 \rangle}

\newcommand{\singsupp}{\mathrm{sing supp}}
\newcommand{\supp}{\mathrm{supp}}

\newfont{\bigmath}{cmr12 at 13pt}

\newfont{\grecomath}{cmmi12 at 15pt}

\newfont{\bl}{msbm10 scaled \magstep2}

\newcommand{\beq}{\begin{equation}}
\newcommand{\eeq}{\end{equation}}
\newcommand{\notmid}{\mid\kern-0.5em\not\kern0.5em}

\newcommand{\eps}{\varepsilon}

\newcommand{\Om}{\Omega}

\newcommand{\val}{\mathrm{v}} 
\newcommand{\esp}{\mathrm{e}}
\newcommand{\M}{\mathcal{M}}

\newcommand{\mP}{\mathcal{P}}

\newcommand{\mO}{\mathcal{O}}

\newcommand{\Gtii}{\Gt^{\infty,\infty}}

\begin{document}
\title{\bf Topological structures in Colombeau algebras: investigation of the duals of $\Gc(\Om)$, $\G(\Om)$ and $\GS(\R^n)$} 
\author{Claudia Garetto \footnote{Current address: Institut f\"ur Technische Mathematik, Geometrie und Bauniformatik,
Universit\"at Innsbruck, e-mail:\,\texttt{claudia@mat1.uibk.ac.at}}\\
Dipartimento di Matematica\\ Universit\`a di Torino, Italia\\
\texttt{garettoc@dm.unito.it}\\
} 
\date{ }
\maketitle

\begin{abstract} We study the topological duals of the Colombeau algebras $\Gc(\Om)$, $\G(\Om)$ and $\GS(\R^n)$, discussing some continuous embeddings and the properties of generalized delta functionals.
\end{abstract}



\setcounter{section}{-1}
\section{Introduction}
This paper is the second part of an investigation into the topological aspects of Colombeau algebras and duality theory started in \cite{Garetto:04b}. While \cite{Garetto:04b} develops the theory of locally convex topological $\wt{\C}$-modules $\G$ and topological duals $\L(\G,\wt{\C})$ in full generality, here we focus on some specific examples: the well-known Colombeau algebras $\Gc(\Om)$, $\G(\Om)$, $\GS(\R^n)$ \cite{Colombeau:85, Colombeau:92, Garetto:04, GGO:03, GKOS:01, NPS:98, O:92} and the corresponding duals $\L(\Gc(\Om),\wt{\C})$, $\L(\G(\Om),\wt{\C})$, $\L(\GS(\R^n),\wt{\C})$. We recall that a locally convex topological $\wt{\C}$-module is a $\wt{\C}$-module $\G$ endowed with a topology $\tau$ which respects the $\wt{\C}$-module structure and has a base of convex neighborhoods of the origin. Here the convexity has to be interpreted as a ``$\wt{\C}$-convexity'' as introduced in \cite[Definition 1.1]{Garetto:04b}. In analogy to locally convex topological vector spaces $\tau$ may be determined by a family of suitable ``$\wt{\C}$-seminorms'' called ultra-pseudo-seminorms whose definition and properties are collected in \cite[Definition 1.8, Theorem 1.10]{Garetto:04b}. These basic notions of topology allow to deal with continuity of $\wt{\C}$-linear maps and to define the topological dual of $\G$ as the space $\L(\G,\wt{\C})$ of all $\wt{\C}$-linear and continuous functionals $T:\G\to\wt{\C}$.

Many results of this paper are obtained as applications of theorems and propositions in \cite{Garetto:04b} while others, more related to Colombeau theory's issues as characterizations of ideals and estimates of representatives, are provided independently. We now describe the contents of the sections in more detail.

Inspired by distribution theory, Section 1 discusses some examples of continuous and $\wt{\C}$-linear functionals as  elements of the topological duals of $\Gc(\Om)$, $\G(\Om)$ and $\GS(\R^n)$. Among them particular attention is given to the direct embedding of classical distributions into $\L(\Gc(\Om),\wt{\C})$ via the map $\D'(\Om)\to\L(\Gc(\Om),\wt{\C}):w\to (u\to [(w(u_\eps))_\eps])$ and to a generalization of the classical Dirac measure, the generalized delta functional $\delta_{\wt{x}}$, whose construction is based on the theory of point values for Colombeau functions \cite{GKOS:01, OK:99}. Since $\Om\to\L(\Gc(\Om),\wt{\C})$ is a sheaf we can talk about the support of a map in $\L(\Gc(\Om),\wt{\C})$ and we manage to characterize the elements of $\L(\G(\Om),\wt{\C})$ as the maps in $\L(\Gc(\Om),\wt{\C})$ having compact support. Under this point of view the analogies between $\L(\G(\Om),\wt{\C})$, $\L(\Gc(\Om),\wt{\C})$ and $\E'(\Om)$, $\D'(\Om)$ respectively are evident. However, relevant differences between the distributional framework and the ``dual-Colombeau'' framework studied here still exist. This fact is pointed out by the computation of the support of $\delta_{\wt{x}}$ and by some remarks on the functionals in $\L(\Gc(\Om),\wt{\C})$ with support $\{0\}$. 

Section 2 is devoted to the properties of the generalized delta functional $\delta_{\wt{x}}$. The main theorem concerns the restriction of $\delta_{\wt{x}}$ to regular generalized functions, the elements of the algebras $\Gcinf(\Om)$, $\Ginf(\Om)$, $\GSinf(\R^n)$ as defined in \cite{Garetto:04, GGO:03, O:92}). It yields that when $\wt{x}$ is a generalized point with compact support (that is, $\wt{x}\in\wt{\Om}_{\rm{c}}$ in the notation of \cite{GKOS:01}), $\delta_{\wt{x}}$ restricted to $\Ginf(\Om)$ is an integral operator of the form $\int_\Om v(y)\cdot\, dy$, where $v$ belongs to $\Gc(\Om)$. It follows that every $u$ in $\Ginf(\Om)$ can be locally represented by the integral $\int_\Om v(x,y)u(y)\, dy$ where $v$ is a suitable generalized function on $\Om\times\Om$. Moreover, an adaptation of these results, involving $\wt{x}\in\wt{\Om}$, $u\in\Gcinf(\Om)$ and $\wt{x}\in\wt{\R^n}$, $u\in\GSinf(\R^n)$, is given and examples show that what we obtained cannot be extended to generalized functions which are not regular.

The concluding Section 3 deals with the continuous embedding of Colombeau algebras into their topological duals. In detail, the chains of inclusions
\[
\Ginf(\Om)\subseteq\G(\Om)\subseteq\L(\Gc(\Om),\wt{\C})
\]
and
\[
\Gc^\infty(\Om)\subseteq\Gc(\Om)\subseteq\L(\G(\Om),\wt{\C})
\]
hold on any open subset $\Om$ of $\R^n$ while algebras of generalized functions satisfying global estimates on $\R^n$ as $\GSinf(\R^n)$, $\GS(\R^n)$ and $\Gt(\R^n)$ \cite{GKOS:01} are naturally contained in $\L(\GS(\R^n),\wt{\C})$ in the following way:
\[
\GSinf(\R^n)\subseteq \GS(\R^n)\subseteq \Gt(\R^n)\subseteq \L(\GS(\R^n),\wt{\C}).
\]
Concerning the proof of this last statement, the Colombeau algebra $\Gpp(\R^n)$ \cite{BO:92} and a useful characterization of the ideals $\NS(\R^n)$ of $\GS(\R^n)$ and $\Npp(\R^n)$ of $\Gpp(\R^n)$ by estimates involving the 0-th derivatives are needed. In addition to such topics we discuss some density issues and we compare the straightforward embedding of classical distributions into $\L(\Gc(\Om),\wt{\C})$ mentioned in Section 1 with the ``Colombeau-embedding''  $\D'(\Om)\subseteq\G(\Om)\subseteq\L(\Gc(\Om),\wt{\C})$. Finally the classical regularization of the Dirac measure via convolution with a mollifier gives the basic idea for a regularization of $\delta_{\wt{x}}\in\L(\GS(\R^n),\wt{\C})$ via generalized functions in $\GS(\R^n)$.

The results obtained in this paper and in particular the chains of inclusions of Colombeau algebras and their duals set up a systematic framework where the regularity theory for pseudodifferential operators with generalized symbols \cite{Garetto:04, GarH:04, GGO:03, GH:03, HO:03} can be settled and developed. The dual $\L(\Gc(\Om\times\Om),\wt{\C})$ was already used in \cite{GGO:03} for defining the kernel of a pseudodifferential operator acting on generalized functions. The dual approach involving $\GSinf(\R^n)$, $\GS(\R^n)$, $\Gt(\R^n)$ and $\L(\GS(\R^n),\wt{\C})$ to pseudodifferential operators with global generalized symbols has been considered in \cite[Chapter 5]{Garetto:04th} for the first time, leading to results of $\GSinf$-regularity which improve the ones obtained in \cite{Garetto:04}.
 
\section{Topological duals of the Colombeau algebras $\Gc(\Om)$, $\G(\Om)$ and $\GS(\R^n)$: basic properties and examples}
In this section we deal with the topological duals of the well-known Colombeau algebras $\Gc(\Om)$, $\G(\Om)$, $\GS(\R^n)$, presenting some basic properties and relevant examples.
Our theoretical background consists of the results on locally convex topological $\wt{\C}$-modules and duality theory presented in the first paper on topological structures in Colombeau algebras \cite{Garetto:04b}. We also refer to \cite{Garetto:04b} for all the notations used in the following.

We recall that for the topologies introduced in \cite[Examples 3.6, 3.7]{Garetto:04b} $\G(\Om)$ and $\GS(\R^n)$ are Fr\'echet $\wt{\C}$-modules while $\Gc(\Om)$ is the strict inductive limit of a sequence of Fr\'echet $\wt{\C}$-modules. Hence by Propositions 2.9 and 2.10 in \cite{Garetto:04b} we know that $\Gc(\Om)$, $\G(\Om)$ and $\GS(\R^n)$ are bornological locally convex topological $\wt{\C}$-modules and that the spaces $\L(\G_c(\Om),\wt{\C})$, $\L(\G(\Om),\wt{\C})$ and $\L(\GS(\R^n),\wt{\C})$ are complete, when endowed with the strong topology $\beta$ or the topology $\beta_b$ of uniform convergence on bounded subsets. Note that by \cite[Propositions 2.14, 2.15]{Garetto:04b} $\Gc(\Om)$, $\G(\Om)$ and $\GS(\R^n)$ are also barrelled. Consequently their topological duals equipped with the weak topology $\sigma$ are quasi-complete. 

Our investigation into the duals of $\Gc(\Om)$, $\G(\Om)$ and $\GS(\R^n)$ begins with some structural issues.

It is clear that if $\Omega'\subseteq\Omega$ then $\Gc(\Om')$ is continuously embedded into $\Gc(\Om)$. This means that every map $T\in\L(\Gc(\Om),\wt{\C})$ can be restricted to $\Om'$ letting 
\beq
\label{restriction}
T_{\vert_{\Om'}}:\Gc(\Om')\to\wt{\C}:u\to Tu,
\eeq
and that $T_{\vert_{\Om'}}$ belongs to $\L(\Gc(\Om'),\wt{\C})$. We are ready now to state the following theorem whose proof is easily obtained by adapting the corresponding result concerning $\D'(\Om)$ to the Colombeau context as in Theorem 1.2.4 \cite{GKOS:01}. Note that Theorem \ref{theorem_sheaf} also follows from abstract arguments in sheaf theory, see e.g. \cite[(2.2.4)]{KaSha:90}. 
\begin{theorem}
\label{theorem_sheaf}
$\Omega\to\L(\Gc(\Om),\wt{\C})$ is a sheaf.
\end{theorem}
It follows that the \emph{support} of $T\in\L(\Gc(\Om),\wt{\C})$ (denoted by $\supp\, T$) can be defined as the complement of the largest open set in $\Om$ in which $T$ is $0$ and that $x\in\Om\setminus\supp\, T$ if and only if there exists an open neighborhood $V\subseteq\Om$ of $x$ such that $T_{\vert_V}=0$ in $\L(\Gc(V),\wt{\C})$. 

In analogy with distribution theory we use the notion of support in order to compare $\L(\Gc(\Om),\wt{\C})$ and $\L(\G(\Om),\wt{\C})$. First of all every $T\in\L(\G(\Om),\wt{\C})$ restricted to $\Gc(\Om)$ is an element of $\L(\Gc(\Om),\wt{\C})$ since the inclusion map $\iota:\Gc(\Om)\to\G(\Om)$ is continuous. Note that $Im(\iota)$ is dense in $\G(\Om)$. In fact for every $u_0\in\G(\Om)$, every seminorm $p_{K,m}(f)=\sup_{x\in K,|\alpha|\le m}|\partial^\alpha f(x)|$ on $\E(\Om)$ and every neighborhood $U:=\{u\in\G(\Om):\ \mP_{K,m}(u-u_0)\le\eta\}$ of $u_0$ in $\G(\Om)$, where $\mP_{K,m}$ is the corresponding ultra-pseudo-seminorm \cite{Garetto:04b}, we can always find some $u\in\Gc(\Om)$ contained in $U$. It suffices to take $\chi\in\Cinf_{\rm{c}}(\Om)$ identically $1$ in a neighborhood of $K$ and observing that $\chi u_0\in\Gc(\Om)$ with $\mP_{K,m}(\chi u_0-u_0)=0$. As a consequence the map
\beq
\label{map_revised}
\L(\G(\Om),\wt{\C})\to\L(\Gc(\Om),\wt{\C}):T\to T_{\vert_{\Gc(\Om)}}
\eeq
is injective.
\begin{theorem}
\label{theorem_compact_support}
A map $T\in\L(\Gc(\Om),\wt{\C})$ belongs to $\L(\G(\Om),\wt{\C})$ if and only if $\supp\, T$ is a compact subset of $\Om$.
\end{theorem}
\begin{proof}
Let $\supp\, T$ be a compact subset of $\Om$ and let $\chi\in\Cinf_{\rm{c}}(\Om)$ identically $1$ in a neighborhood of $\supp\, T$. For all $u\in\G(\Om)$ we have that $\chi u\in\Gc(\Om)$ and 
\[
T':\G(\Om)\to\wt{\C}:u\to T(\chi u)
\]
is a $\wt{\C}$-linear map on $\G(\Om)$. Moreover, the restriction of $T'$ to $\Gc(\Om)$ coincides with $T$ since $(1-\chi)u\in\Gc(\Om\setminus\supp\, T)$ and $T'(u)-T(u)=T((1-\chi)u)=0$. Note that for all $u\in\G(\Om)$, $\chi u$ is a generalized function in $\G_K(\Om)$ with $K=\supp\, \chi$ and by continuity of $T$ the estimate 
\beq
\label{est_continuity_1}
\vert T'(u)\vert_\esp\le C\mP_{\G_K(\Om),m}(\chi u),
\eeq
with $\mP_{\G_K(\Om),m}$ introduced in \cite[Example 3.7]{Garetto:04b}, is valid on $\G(\Om)$. By the Leibniz rule $\mP_{\G_K(\Om),m}(\chi\cdot)$ can be bounded by some ultra-pseudo-seminorm which determines the sharp topology on $\G(\Om)$. Combined with \eqref{est_continuity_1} this guarantees the continuity of $T'$ as a $\wt{\C}$-linear map from $\G(\Om)$ into $\wt{\C}$.\\
Assume now that $T\in\L(\G(\Om),\wt{\C})$ but that its support is not a compact subset of $\Om$. Let $K_0\subset K_1\subset K_2\subset...$ be an exhausting sequence of compact sets of $\Om$. Then for all $n\in\N$, $\supp\, T\cap (\Om\setminus K_n)\neq \emptyset$. In other words there exists a sequence $(u_n)_n\subseteq \Gc(\Om)$ such that $\supp\, u_n\subseteq \Omega\setminus K_n$ and $Tu_n\neq 0$. Denoting the valuation of $Tu_n\in\wt{\C}$ by $a_n$, the generalized function $v_n=[(\eps^{-a_n})_\eps]u_n$ has support contained in $\Om\setminus K_n$ and $\vert Tv_n\vert_\esp=\esp^{a_n}\vert Tu_n\vert_\esp= 1$. The sequence $(v_n)_n$ converges to $0$ in $\G(\Om)$ since for all $K\Subset\Om$ there exists $n_0\in\N$ such that $K\subseteq K_{n_0}$ and from $K_{n_0}\subseteq \Om\setminus \supp\, v_n$ for all $n\ge n_0$ we obtain that 
\[
\sup_{x\in K,|\alpha|\le m}|\partial^\alpha v_{n,\eps}(x)|\le \sup_{x\in K_{n_0},|\alpha|\le m}|\partial^\alpha v_{n,\eps}(x)|=O(\eps^q),
\]
where $q$ is any natural number and $n\ge n_0$. This contradicts the continuity of $T$ on $\G(\Om)$.
\end{proof}
From Theorem \ref{theorem_compact_support} it follows that the restriction map in \eqref{map_revised}
is a bijection between $\L(\G(\Om),\wt{\C})$ and the space of all functionals in $\L(\Gc(\Om),\wt{\C})$ with compact support.

Distribution theory and point value theory for generalized functions are the main tools in providing the following examples of $\wt{\C}$-linear and continuous functionals on $\Gc(\Om)$, $\G(\Om)$ and $\GS(\R^n)$. The natural embeddings of Colombeau algebras into their topological duals are postponed to Section 3.
\begin{example}
\label{example_distributions}
 
\bf{Classical distributions}\rm

Every distribution $w\in\D'(\Om)$ provides an element of $\L(\Gc(\Om),\wt{\C})$. In fact recalling that the injection $\iota_\D:\Gc(\Om)\to\G_{\D(\Om)}: u\to (u_\eps)_\eps+\Neg_{\D(\Om)}$ introduced in \cite[Remark 3.8]{Garetto:04b} is continuous for the strict inductive limit topology on $\Gc(\Om)$ and the sharp topology on $\G_{\D(\Om)}$, we have by \cite[(ii)Remark 3.14]{Garetto:04b} that $$w\circ\iota_\D:\G_c(\Om)\to\wt{\C}:w\to[(w(u_\eps))_\eps],$$ where $(u_\eps)_\eps$ is any representative of $u$ in some $\M_{\D_{K'}(\Om)}$, is $\wt{\C}$-linear and continuous. In the same way we get an element of $\L(\G(\Om),\wt{\C})$ and an element of $\L(\GS(\R^n),\wt{\C})$ if $w$ belongs to $\E'(\Om)$ and $\S'(\R^n)$ respectively.
\end{example}
In the following examples we make use of the main definitions and results of generalized point value theory. The reader should refer to \cite[Subsection 1.2.4]{GKOS:01} for the needed notations and explanations.
\begin{example}
\label{example_delta}

\bf{The generalized delta functional $\delta_{\wt{x}}$}\rm

Let $\wt{x}$ be a generalized point in $\wt{\Om}_{\rm{c}}$. We can define a $\wt{\C}$-linear map $\delta_{\wt{x}}:\G(\Om)\to\wt{\C}$ associating with each $u\in\G(\Om)$ its point value $u(\wt{x})$ at $\wt{x}$. It is clear that $\delta_{\wt{x}}$ belongs to $\L(\G(\Om),\wt{\C})$ since 
\[
\vert\delta_{\wt{x}}(u)\vert_\esp\le \esp^{-\val(\sup_{x\in K}|u_\eps(x)|)}
\]
where $(x_\eps)_\eps$ is a representative of $\wt{x}$ contained in a compact set $K$ for $\eps$ small enough.

Since every $u\in\Gc(\Om)$ is a generalized function in some $\G_{\D_{K'}(\Om)}$ it is meaningful to define the point value $u(\wt{x})$ even when $\wt{x}$ is in $\wt{\Om}\setminus\wt{\Om}_{\rm{c}}$. Indeed for $(u_\eps)_\eps\in\M_{\D_{K'}(\Om)}$ and $(x_\eps)_\eps\in\Om_M$ the estimate $|u_\eps(x_\eps)|\le\sup_{y\in\Omega}|u_\eps(y)|$ tells us that $u(\wt{x}):=[(u_\eps(x_\eps))_\eps]$ is well-defined for every choice of representative of $u$ of this type. Therefore, $\delta_{\wt{x}}$ is a $\wt{\C}$-linear map from $\Gc(\Om)$ into $\wt{\C}$ for each $\wt{x}\in\wt{\Om}$ and by $\vert\delta_{\wt{x}}(u)\vert_\esp\le\mP_{K,0}(u)$ on $\G_K(\Om)$ it is continuous.
\end{example}
We conclude this section with some further considerations and examples involving the generalized delta functional. They will point out some differences between the distributional delta and the generalized delta concerning topics such as support and representation of continuous functionals.
Some interesting phenomena occur when we compute the support of $\delta_{\wt{x}}\in\L(\Gc(\Om),\wt{\C})$.
\begin{example}
\bf{The support of $\delta_{\wt{x}}$}\rm

Let $\wt{x}=[(x_\eps)_\eps]\in\wt{\Om}$. Defining the support $\supp\, \wt{x}$ of $\wt{x}$ as the complement of the set
\beq
\label{supp_wtx}
\left\{ x_0\in\Om:\ \begin{array}{cc}\exists V(x_0)\ \text{neighborhood of $x_0 $}\ \ \exists \eta\in(0,1]\ \ \forall\eps\in(0,\eta]\\[0.2cm]
x_\eps\not\in V(x_0)
\end{array}
\right\},  
\eeq
in $\Om$, we want to prove that $\supp\, \delta_{\wt{x}}=\supp\, \wt{x}$.
The reader can easily check that $\Om\setminus\supp\, \wt{x}$ is well-defined since it does not depend on the representative $(x_\eps)_\eps$ of $\wt{x}$.
% Indeed if $(x_\eps)_\eps$ satisfies \eqref{supp_wtx} from the fact that $x_\eps\not\in V(x_0)$ for all $\eps\in(0,\eta]$ it follows that $|x_\eps-x_0|>r$ for some $r>0$. If $|x_\eps-x'_\eps|=O(\eps^q)$ for all $q\in\N$ choosing possibly a smaller $\eta'$ we have the estimate
%\[
%|x'_\eps-x_0|=|x'_\eps-x_\eps+x_\eps-x_0|\ge |x_\eps-x_0|-|x'_\eps-x_\eps|>r-c\eps>\frac{r}{2}.
%\]
%on $(0,\eta']$. This means that $x'_\eps$ does not belong to $V'(x_0):=\{x\in\Om:\ |x-x_0|\le r/2\}$ when $\eps$ is varying in the interval $(0,\eta']$. 
If $x_0\in\Om\setminus\supp\ {\wt{x}}$ then there exists a neighborhood $V(x_0)$ satisfying \eqref{supp_wtx} and we can take a smaller $V'(x_0)$ such that $V'(x_0)\subseteq\overline{V'(x_0)}\Subset V(x_0)$. Any generalized function $u\in\Gc(V'(x_0))$ has a representative $(u_\eps)_\eps$ with $\supp\, u_\eps\subseteq V(x_0)$ for all $\eps\in(0,1]$. Hence by \eqref{supp_wtx} $\delta_{\wt{x}}(u)=0$ for all $u\in\Gc(V'(x_0))$ and therefore, $x_0\in\Om\setminus\supp\, \delta_{\wt{x}}$. Conversely if $x_0\in\Om\setminus\supp\,\delta_{\wt{x}}$ then there exists a neighborhood $V(x_0)$ such that $\delta_{\wt{x}}(u)=0$ for all $u\in\Gc(V(x_0))$. In particular, for all $\psi\in\Cinf_{\rm{c}}(V(x_0))$, representatives $(x_\eps)_\eps$ of $\wt{x}$ and $q\in\N$,
\beq
\label{neg_psi}
|\psi(x_\eps)|=O(\eps^q)\ \ \text{as}\ \eps\to 0.
\eeq
This implies that $x_0\in\Om\setminus\supp\ \wt{x}$ since taking $V'(x_0)\subseteq\overline{V'(x_0)}\Subset V(x_0)$ and $\psi\in\Cinf_{\rm{c}}(V(x_0))$ identically $1$ on $V'(x_0)$, \eqref{neg_psi} can hold only if $x_\eps$ does not belong to $V'(x_0)$ for $\eps$ small.
%In fact if we assume that $x_0\in\supp\ \wt{x}$ then we can find a sequence $(\eps_n)_n$ tending to $0$ such that $(x_{\eps_n})_n$ is contained in some $V'(x_0)\subseteq\overline{V'(x_0)}\Subset V(x_0)$. Hence for $\psi\in\Cinf_{\rm{c}}(V(x_0))$ identically $1$ on $V'(x_0)$, we have that $\psi(x_{\eps_n})=1$ for all $n$ in contradiction with \eqref{neg_psi}.

When $\wt{x}\in\wt{\Om}_{\rm{c}}$ it is clear that $\supp\ \delta_{\wt{x}}$ is a compact subset of $\Om$ since $\supp\ \wt{x}$ is contained in some $K\Subset\Om$. 

Note that by \eqref{supp_wtx} the support of $\wt{x}\in\wt{\Om}$ is the set of all accumulation points of a representing net $(x_\eps)_\eps$. As a consequence and differently from the distributional case, $\delta_{\wt{x}}$ may have support empty or unbounded. Consider in $\R$ the example $\delta_{[(\eps^{-1}x)_\eps]}$ where $x\neq 0$ is some point of $\R$. Since $\eps^{-1}|x|\to +\infty$ we have that $\supp(\delta_{[(\eps^{-1}x)_\eps]})=\emptyset$. This shows that $\delta_{[(\eps^{-1}x)_\eps]}$ is identically null on $\Gc(\R)$.
Finally let $(x_n)_n$ be a sequence in $\R^+$ with $\N$ as set of accumulation points and $x_n\le n+1$ for all $n\in\N$. Defining $x_\eps=x_{n}$ for $\eps\in(1/(n+2),1/(n+1)]$, $(x_\eps)_\eps$ is moderate since $x_\eps=x_{n}\le n+1\le\eps^{-1}$ on each interval $(1/(n+2),1/(n+1)]$. Thus for $\wt{x}:=[(x_\eps)_\eps]\in\wt{\R}$ we conclude that $\supp\, \delta_{\wt{x}}=\N$. 
\end{example}
\begin{remark}
\label{remark_no_distrib}
Another difference to distribution theory consists in the existence of elements of $\L(\G(\Om),\wt{\C})$ having support $\{0\}$ which are not $\wt{\C}$-linear combinations of $\delta_0$ and its derivatives. As an example take $[(\eps)_\eps]\in\wt{\R}$ and $\delta_{[(\eps)_\eps]}\in\L(\G(\R),\wt{\C})$. By \eqref{supp_wtx} $\supp\, \delta_{[(\eps)_\eps]}=\{0\}$ but for any $m\in\N$ and any choice of $c_i\in\wt{\C}$ the equality
\beq
\label{equality_no_distrib}
\delta_{[(\eps)_\eps]}=\sum_{i=0}^m c_i\delta_{0}^{(i)}
\eeq
does not hold. In fact if we assume that $\delta_{[(\eps)_\eps]}$ is equal to $\sum_{i=0}^m c_i\delta_{0}^{(i)}$ in $\L(\G(\R),\wt{\C})$ then taking $\psi\in\Cinfc(\R)$ identically $1$ in a neighborhood of the origin we obtain $\delta_{[(\eps)_\eps]}(\psi)=c_0\psi(0)$ which implies $1=[(\psi(\eps))_\eps]=c_0$. Moreover, the application of \eqref{equality_no_distrib} to each function $x^h\psi(x)$ with $h\in\N$, $1\le h\le m$, yields $c_i=[((-\eps)^i/i!)_\eps]$ for $1\le i\le m$. Hence $\delta_{[(\eps)_\eps]}$ has to coincide with $\sum_{i=0}^{m}[((-\eps)^i/i!)_\eps]\delta^{(i)}_0$. But this is impossible since for $x^{m+1}\psi(x)$ we have that $\delta_{[(\eps)_\eps]}(x^{m+1}\psi(x))=[(\eps^{m+1})_\eps]\neq 0$ while $\sum_{i=0}^{m}[((-\eps)^i/i!)_\eps]\delta^{(i)}_0(x^{m+1}\psi(x))=0$ in $\wt{\C}$.

Note that $\delta_{[(\eps)_\eps]}$ can be written as the series $\sum_{i=0}^{\infty}[((-\eps)^i/i!)_\eps]\delta^{(i)}_0$ when it acts on $\Ginf(\R)$, in the sense that for all $u\in\Ginf(\R)$
\beq
\label{convergence_formula}
u([(\eps)_\eps])=\sum_{i=0}^\infty [(\eps^i/i!)_\eps]u^{(i)}(0)\qquad \text{in}\ \wt{\C}.
\eeq
In fact if $(u_\eps)_\eps\in\E^\infty_M(\R)$, \eqref{convergence_formula} follows from 
\[
\biggl|u_\eps(\eps)-\sum_{i=0}^q\frac{\eps^i}{i!}u_\eps^{(i)}(0)\biggr|=\biggl|\frac{\eps^{q+1}}{(q+1)!}u_\eps^{(q+1)}(\theta\eps)\biggr|=O(\eps^{q+1-N}),
\]
where $N$ does not depend on $q$.

Take now $f\in\Cinf(\R)$ real valued with $\sum_{i=0}^\infty f^{(i)}(0)/i!$ divergent. The generalized function $u:=[(f(\cdot/\eps))_\eps]$ belongs to $\G(\R)\setminus \Ginf(\R)$ and by construction the series $\sum_{i=0}^\infty[((-\eps)^i/i!)_\eps]\delta_0^{(i)}(u)=\sum_{i=0}^\infty [(f^{i}(0)/i!)_\eps]$ is not convergent in $\wt{\C}$. This proves that \eqref{convergence_formula} can not be strengthened to hold in $\G(\R)$.
\end{remark}
We recall that Corollary 2.18 in \cite{Garetto:04b} valid for barrelled locally convex topological $\wt{\C}$-modules, says in particular that when we have a sequence $(T_n)_n$ in $\L(\GS(\R^n),\wt{\C})$ with $(T_n(u))_n$ convergent to some $T(u)$ in $\wt{\C}$ for all $u$, then the limit $T$ belongs to $\L(\GS(\R^n),\wt{\C})$. Analogously, the same holds for sequences in $\L(\Gc(\Om),\wt{\C})$ and $\L(\G(\Om),\wt{\C})$. We can use this convergence result for constructing the following example of a series of generalized delta functionals. Note that when $\wt{x}\in\wt{\R^n}$ the point value $u(\wt{x}):=[(u_\eps(x_\eps))_\eps]$ of $u\in\GS(\R^n)$ is a well-defined element of $\wt{\C}$ and the delta functional $\delta_{\wt{x}}:u\to u(\wt{x})$ is a map in $\L(\GS(\R^n),\wt{\C})$.
\begin{example}
\label{example_sequence}
Consider the sequence of generalized points $\wt{x}_n:=[(\eps^{-n})_\eps]$ in $\wt{\R}$. We know that $\delta_{\wt{x}_n}\in\L(\GS(\R),\wt{\C})$ and the series $\sum_{n=0}^\infty u(\wt{x}_n)$ is convergent in $\wt{\C}$. This is due to the fact that if $\sup_{x\in\R}(1+|x|)|u_\eps(x)|=O(\eps^{-N})$ as $\eps$ tends to $0$ then
\[
\biggl|\sum_{n=q+N}^m u_\eps(\eps^{-n})\biggr|\le c\eps^{-N}\sum_{n=q+N}^m \eps^n\le c'\eps^q
\]
for $\eps$ small enough, or in other words $\val_{\wt{\C}}(\sum_{n=q+N}^m\delta_{\wt{x}_n}(u))\ge q$. Therefore, by the corollary quoted above the map $u\to\sum_{n=0}^\infty\delta_{\wt{x}_n}(u)$ belongs to $\L(\GS(\R),\wt{\C})$.
\end{example} 
\section{Integral representation of $\delta_{\wt{x}}$ over regular Co\-lom\-be\-au functions}
In the previous section we presented the generalized delta functional $\delta_{\wt{x}}$ as a typical example of a $\wt{\C}$-linear map in $\L(\Gc(\Om),\wt{\C})$. We continue to investigate the properties of $\delta_{\wt{x}}$ by considering the action on regular Colombeau generalized functions. Since we will work with integral operators we recall that a generalized function $u\in\G(\Om)$ can be integrated over a compact set $K$ of $\Om$ defining $\int_K u(y)\, dy$ as the complex generalized number  with representative $(\int_K u_\eps(y)\, dy)_\eps$. When $u\in\Gc(\Om)$ we set $\int_\Om u(y)\, dy :=\int_{K}u(y)\, dy$ where $K$ is any compact set containing $\supp\, u$ in its interior. We are now ready to prove the following theorem.
\begin{theorem}
\label{delta_theorem}
For all $\wt{x}\in\wt{\Om}_{\rm{c}}$ there exists $v\in\Gc(\Om)$ such that for all $u\in\Ginf(\Om)$
\beq
\label{int_representation_formula}
u(\wt{x})=\int_\Om v(y)u(y)\, dy.
\eeq
\end{theorem}
\begin{proof}
We fix $\wt{x}\in\wt{\Om}_{\rm{c}}$. From the definition of a compactly supported generalized point there exists a compact set $K\Subset\Om$ such that every representative $(x_\eps)_\eps$ of $\wt{x}$ is contained in $K$ for $\eps$ small enough. Let $\psi$ be a function in $\Cinfc(\Om)$ identically $1$ on a neighborhood of $K$ and $\varphi$ a mollifier in $\S(\R^n)$ with $\int\varphi =1$ and $\int y^\alpha \varphi =0$ for all $\alpha\neq 0$ in $\N^n$. We prove that 
\beq
\label{def_function_v}
v:=(\psi(\cdot)\varphi_\eps(x_\eps-\cdot))_\eps+\Neg(\Om)
\eeq
is a well-defined generalized function in $\Gc(\Om)$. In fact for $(x_\eps)_\eps\in\Om_M$ we have that $(\psi(y)\varphi_\eps(x_\eps-y))_\eps\in\EM(\Om)$ and by a Taylor's formula argument as in \cite[Proposition 1.2.45]{GKOS:01} if $(x_\eps)_\eps\sim(x'_\eps)_\eps$ then $(\psi(y)(\varphi_\eps(x_\eps-y)-\varphi_\eps(x'_\eps-y)))_\eps$ belongs to $\Neg(\Om)$. Thus $v\in\G(\Om)$ with $\supp\, v\subseteq \supp\, \psi\Subset \Om$. 
%we can write
%\[
%\varphi_\eps(x_\eps-y)-\varphi_\eps(x'_\eps-y)=\nabla\varphi_\eps(x'_\eps-y+\theta(x_\eps-x'_\eps))(x_\eps-x'_\eps),
%\]
%where $\theta\in[0,1]$. As a consequence for any $q\in\N$ there exists $\eta\in(0,1]$ such that for all $\eps\in(0,\eta]$ and for all $y\in\R^n$
%\beq
%\label{estimate_function_v_1}
%|\varphi_\eps(x_\eps-y)-\varphi_\eps(x'_\eps-y)|\le c\eps^{-n-1}|x_\eps-x'_\eps|\le \eps^{-n-1+q}.
%\eeq
%Since $(\varphi_\eps(x_\eps-y)-\varphi_\eps(x'_\eps-y))_\eps\in\EM(\R^n)$, \eqref{estimate_function_v_1} means that $(\varphi_\eps(x_\eps-y)-\varphi_\eps(x'_\eps-y))_\eps\in\Neg(\R^n)$ which implies $(\psi(y)\varphi_\eps(x_\eps-y)-\psi(y)\varphi_\eps(x'_\eps-y))_\eps\in\Neg(\Om)$. This results allows us to say that $v\in\G(\Om)$ with $\supp\, v\subseteq \supp\, \psi\Subset\Om$ by construction.
We complete the proof by showing that
\[
\biggl(u_\eps(x_\eps)-\int_\Om u_\eps(y)\psi(y)\varphi_\eps(x_\eps-y)\, dy\biggr)_\eps\in\Neg.
\]
First of all the net $u'_\eps(y)=u_\eps(y)\psi(y)$ defines an element of $\E^\infty_M(\R^n)$ with $\supp\, u'_\eps\subseteq\supp\, \psi$ for all $\eps\in(0,1]$ and the equality
\beq
\label{equality_v_1}
\int_\Om u_\eps(y)\psi(y)\varphi(x_\eps-y)\, dy =\int_{\R^n}u'_\eps(x_\eps-\eps y)\varphi(y)\, dy
\eeq
holds. Taylor expansion of $u'_\eps$ at $x_\eps$ and the properties of $\varphi$ yield the estimate
\beq
\label{estimate_function_v_2}
\biggl|\int_{\R^n}u'_\eps(x_\eps-\eps y)\varphi(y)\, dy -u'_\eps(x_\eps)\biggr|\le c\eps^{-N+q+1},
\eeq
for $\eps$ small enough and arbitrary $q\in\N$, where $N$ is independent of $q$ by the $\E^\infty_M$-property of $(u_\eps)_\eps$. Note that again for $\eps$ small, the point $x_\eps$ belongs to $K$, $\psi(x_\eps)=1$ and then $u_\eps(x_\eps)=u'_\eps(x_\eps)$. Combining \eqref{equality_v_1} with \eqref{estimate_function_v_2} we conclude that for all $q\in\N$ there exists $\eta\in(0,1]$ such that for all $\eps\in(0,\eta]$
\[
\biggl|u_\eps(x_\eps)-\int_\Om u_\eps(y)\psi(y)\varphi_\eps(x_\eps-y)\, dy\biggr|\le \eps^q.
\]
\end{proof}
\begin{remark}
\label{remark_v}
In this remark we collect some comments concerning the generalized function $v$ constructed in the previous theorem.
\begin{trivlist}
\item[$(i)$] The definition of $v$ is independent of the choice of the cut-off function $\psi$. Let $\psi_1$ and $\psi_2$ be two different cut-off functions identically $1$ in a neighborhood of $K\Subset\Om$ and assume that $x_\eps\in K$ for all $\eps\in(0,\eta]$. We know that there exists $\eps_0>0$ such that if $|x_\eps-y|\le\eps_0$ then $\psi_1(y)-\psi_2(y)=0$. This means that it suffices to consider the subset $\{y\in K_1:\ \exists\eps\in(0,\eta]\  |x_\eps-y|>\eps_0\}$ for estimating the net $((\psi_1-\psi_2)(y)\varphi_\eps(x_\eps-y))_\eps$ on a compact set $K_1\subset\Om$ for $\eps\in(0,\eta]$. From the rapidly decreasing behavior of $\varphi$ we have for arbitrary $q\in\N$
\beq
\label{rap_decreasing_varphi}
|\psi_1(y)-\psi_2(y)||\varphi_\eps(x_\eps-y)|\le c\eps^{-n}\lara{\frac{x_\eps-y}{\eps}}^{-q}\le c\eps^{-n}\lara{\frac{\eps_0}{\eps}}^{-q}\le c\eps_0^{-q}\eps^{-n+q}
\eeq
in that region, with $\lara{z}=(1+|z|^2)^{\frac{1}{2}}$. This shows that $(\psi_1(y)\varphi_\eps(x_\eps-y)-\psi_2(y)\varphi_\eps(x_\eps-y))_\eps\in\Neg(\Om)$.
However, it easy to check that the generalized function $v$ may depend on the mollifier $\varphi$.
%\item[$(ii)$] The function $v$ may depend on the choice of the mollifier $\varphi$. In fact taking mollifiers $\varphi_1$ and $\varphi_2$ in $\S(\R^n)$ with $\varphi_1(x_0)\neq\varphi_2(x_0)$ for some point $x_0\neq 0$ of $\Om$, the generalized function
%\[
%w:=\big[\big(\eps^{-n}\psi(y)\big(\varphi_1\big(\frac{x_\eps-y}{\eps}\big)-\varphi_2\big(\frac{x_\eps-y}{\eps}\big)\big)\big)_\eps\big]\in\Gc(\Om)
%\]
%is not $0$. This assertion is clear, since defining the generalized point $\wt{x}'\in\wt{\Om}_{\rm{c}}$ with representative $x_\eps-\eps x_0$ for $\eps$ small enough, we have that
%\[
%w({\wt{x}}')=[(\eps^{-n}\psi(x_\eps-\eps x_0)(\varphi_1-\varphi_2)(x_0))_\eps]=[(\eps^{-n}(\varphi_1-\varphi_2)(x_0))_\eps]\neq 0.
%\]
\item[$(ii)$] In general \eqref{int_representation_formula} does not hold when we take $v$ defined in \eqref{def_function_v} and $u\in\G(\Om)\setminus\Ginf(\Om)$. Let us consider the one-dimensional case for simplicity. We prove that for all mollifiers $\varphi\in\S(\R)$ with $\Vert\varphi\Vert^2_2-\overline{\varphi}(0)\neq 0$ and for all $x\in\R$ there exists $u\in\G(\R)$ such that for all cut-off functions $\psi$ identically $1$ in \ neighborhood of $x$,
\[
u(x)\neq \biggl[\biggl(\int_{\R}u_\eps(y)\psi(y)\varphi_\eps(x-y)\, dy\biggr)_\eps\biggr]\qquad \text{in}\ \wt{\C}.
\]
It is sufficient to define $u$ as the generalized function in $\G(\R)$ with representative $u_\eps(y):=\overline{\varphi}_\eps(x-y)$ and write
\begin{multline*}
\int_{\R}u_\eps(y)\psi(y)\varphi_\eps(x-y)\, dy =\int_{\R}\psi(y)|\varphi_\eps|^2(x-y)\, dy\\ = \int_{\R}|\varphi_\eps|^2(x-y)\, dy + \int_\R(\psi(y)-1)|\varphi_\eps|^2(x-y)\, dy.
\end{multline*}
The same kind of reasoning as used in \eqref{rap_decreasing_varphi} yields $(\int_\R(\psi(y)-1)|\varphi_\eps|^2(x-y)\, dy)_\eps\in\Neg$. Therefore, we conclude from the hypothesis on $\varphi$ that 
\[
\biggl[\biggl(\int_{\R}u_\eps(y)\psi(y)\varphi_\eps(x-y)\, dy\biggr)_\eps\biggr]-u(x)=[(\eps^{-1}\Vert\varphi\Vert_2^2-\eps^{-1}\overline{\varphi}(0))_\eps]\neq 0.
\]
\item[$(iii)$] Finally, dealing with elements of $\Gcinf(\Om)$ we can claim that for all $\wt{x}\in\wt{\Om}$ there exists $v\in\G(\Om)$ such that for all $u\in\Gcinf(\Om)$
\[
u(\wt{x})=\int_\Om v(y)u(y)\, dy.
\]
The proof of Theorem \ref{delta_theorem} shows that this result is valid with $v$ taking the form of $(\varphi_\eps(x_\eps-y)_{\vert_{\Om}})_\eps+\Neg(\Om)$.
\end{trivlist}
\end{remark}  
In the sequel we adopt the notation $v_{\wt{x}}$ for a generalized function $v$ satisfying the statements of Theorem \ref{delta_theorem} or $(iii)$ Remark \ref{remark_v}. (For brevity we omit the dependence on the mollifier $\varphi$).  Since the topology with which we endow $\Gcinf(\Om)$ is finer than the topology induced by $\Gc(\Om)$ on $\Gcinf(\Om)$ we know that the restriction $T_{\vert_{\Gcinf(\Om)}}$ of a map $T\in\L(\Gc(\Om),\wt{\C})$ is an element of $\L(\Gcinf(\Om),\wt{\C})$. In the same way $T_{\vert_{\Ginf(\Om)}}\in\L(\Ginf(\Om),\wt{\C})$ if $T\in\L(\G(\Om),\wt{\C})$. In the particular case of the generalized delta functional $\delta_{\wt{x}}$, Theorem \ref{delta_theorem} and Remark \ref{remark_v} tell us that this restriction is an integral operator of the form $\int_\Om v_{\wt{x}}(y)\cdot \,dy$. This operator is an element of $\L(\Gcinf(\Om),\wt{\C})$ when $\wt{x}\in\wt{\Om}$ and an element of $\L(\Ginf(\Om),\wt{\C})$ when $\wt{x}\in\wt{\Om}_{\rm{c}}$.

Theorem \ref{delta_theorem} leads to a local integral representation for regular Colombeau generalized functions as an interesting consequence.
\begin{theorem}
\label{representation_theorem}
For all open subsets $\Om,\Om'$ of $\R^n$ with $\Om'\subseteq\overline{\Om'}\Subset\Om$ there exists $v\in\G(\Om'\times\Om)$ with $\cup_{x\in\Om'}\supp\, v(x,\cdot)\subseteq\Om$ compact, such that for all $u\in\Ginf(\Om)$
\beq
\label{formula_representation_1}
u_{\vert_{\Om'}}(x)=\int_\Om v(x,y)u(y)\, dy\qquad \text{in}\ \G(\Om').
\eeq
Moreover, there exists $v\in\G(\Om\times\Om)$ such that for all $u\in\Gcinf(\Om)$
\beq
\label{formula_representation_2}
u(x)=\int_\Om v(x,y)u(y)\, dy\qquad \text{in}\ \G(\Om).
\eeq
\end{theorem}
\begin{proof}
The main idea of this proof is to consider $v_x$ where $x\in\Om'$ and to observe that the map $x\to v_x$ defines a generalized function in $\G(\Om'\times\Om)$ fulfilling \eqref{formula_representation_1}. For a cut-off function $\psi$ identically $1$ in a neighborhood of $\overline{\Om'}$ and a mollifier $\varphi$ in $\S(\R^n)$, let us define $v$ as the class of the net $((\psi(y)\varphi_{\eps}(x-y))_{\vert_{\Om'\times\Om}})_\eps$ in $\G(\Om'\times\Om)$. This is possible since the net $((\psi(y)\varphi_{\eps}(x-y))_{\vert_{\Om'\times\Om}})_\eps$ is moderate. Note that $\supp\, v(x,\cdot)\subseteq \supp\, \psi$ for all $x\in\Om'$. By Proposition 2.14 in \cite{GGO:03} the integral $\int_{\Om}v(\cdot,y)u(y)\, dy$ gives a generalized function in $\G(\Om')$. We recall that for all $\wt{x}\in\wt{\Om'}_{\rm{c}}$, the generalized function $v(\wt{x},\cdot):=(v_\eps(x_\eps,y))_\eps+\Neg(\Om)$ belongs to $\Gc(\Om)$ and by construction it coincides with $v_{\wt{x}}$ of Theorem \ref{delta_theorem}. Hence the equality
\[
u_{\vert_{\Om'}}(\wt{x})=\int_{\Om}v_{\wt{x}}(y)u(y)\, dy =\int_{\Om}v(\wt{x},y)u(y)\, dy
\]
is valid for $u\in\Ginf(\Om)$ and proves that $u_{\vert_{\Om'}}$ and $\int_\Om v(\cdot,y)u(y)\, dy$ are the same functions in $\G(\Om')$.\\
If $u\in\Gcinf(\Om)$ then by $(iii)$ Remark \ref{remark_v} we deduce that $v:=[(\varphi_\eps(x-y)_{\vert{\Om\times\Om}})_\eps]\in\G(\Om\times\Om)$ satisfies \eqref{formula_representation_2}.
\end{proof}
We conclude this section with global versions of Theorems \ref{delta_theorem} and \ref{representation_theorem} which involve a suitable subalgebra of regular and tempered generalized functions. Our interest for such a matter is motivated by the relevant role that algebras of generalized functions globally defined on $\R^n$ and the corresponding duals play in some recent papers on generalized pseudodifferential operators and regularity theory \cite{Garetto:04, Garetto:04th}. We will refer again to the following results in Subsection 3.2, once we will have studied the relationships between $\GS(\R^n)$, $\Gt(\R^n)$ and $\L(\GS(\R^n),\wt{\C})$.

Let $\Et[\R^n]:=\OM(\R^n)^{(0,1]}$ and $\Nt(\R^n)$ the ideal of $\Gt(\R^n)$ as defined in \cite[Example 3.9]{Garetto:04b}. We denote by $\Gtii(\R^n)$ the quotient $\Et^{\infty,\infty}(\R^n)/\Nt(\R^n)$ defined by
\begin{multline*}
\Et^{\infty,\infty}(\R^n):=\{(u_\eps)_\eps\in\Et[\R^n]:\  \exists N\in\N\, \forall \alpha\in\N^n\\ \sup_{x\in\R^n}(1+|x|)^{-N}|\partial^\alpha u_\eps(x)|=O(\eps^{-N})\ \text{as}\ \eps\to 0\}.
\end{multline*}
By construction $\Gtii(\R^n)$ is a subalgebra of $\Gt(\R^n)$ which contains $\mO_C(\R^n)$ but not $\mO_M(\R^n)$. We recall that the product of $v\in\Gt(\R^n)$ with $u\in\GS(\R^n)$ is a well-defined element $vu:=(v_\eps u_\eps)_\eps+\NS(\R^n)$ of $\GS(\R^n)$ and that we can integrate a generalized function $w\in\GS(\R^n)$ on $\R^n$; the integral is defined as the complex generalized number $\int_{\R^n}w(y)\, dy :=[(\int_{\R^n}w_\eps(y)\, dy)_\eps]$. Finally when $v\in\Gt(\R^{2n})$ has a representative $(v_\eps)_\eps$ satisfying the condition
\begin{multline}
\label{condvrepr}
\forall\alpha\in\N^n\ \forall k\in\N\ \exists N\in\N\\
\sup_{x\in\R^n}(1+|x|)^{-N}\hskip-8pt\sup_{y\in\R^n, |\beta|\le k}(1+|y|)^k|\partial^\alpha_x\partial^\beta_y v_\eps(x,y)|=O(\eps^{-N}),
\end{multline}
the generalized function $v(x,\cdot):=(v_\eps(x,\cdot))_\eps+\NS(\R^n)$ is a well-defined element of $\GS(\R^n)$ for all $x\in\R^n$, since it does not depend on the choice of the representatives of $v$ fulfilling \eqref{condvrepr}. Then for $u\in\Gt(\R^n)$ the integral $\int_{\R^n}v(x,y)u(y)\, dy$ makes sense for all $x$ and by \eqref{condvrepr} provides a generalized function $\int_{\R^{n}}v(x,y)u(y)\, dy:=(\int_{\R^n}v_\eps(x,y)u_\eps(y)\, dy)_\eps+\Nt(\R^n)$ in $\Gt(\R^n)$. Finally we consider $\Gt(\R^n)$ endowed with the topology determined via $\GtS(\R^n)$ in \cite[Example 3.9]{Garetto:04b} and $\Gtii(\R^n)$ with the topology induced by $\Gt(\R^n)$. Since for $\wt{x}\in\R^n$ and for all $m,N\in\N$ the $\wt{\C}$-linear map $\delta_{\wt{x}}:\G^m_{N,\S}(\R^n)\to\wt{\C}:u\to u(\wt{x})$ is continuous we have that $\delta_{\wt{x}}$ is an element of $\L(\Gt(\R^n),\wt{\C})$.

We now have all the tools for stating and proving Propositions \ref{delta_proposition} and \ref{representation_prop}.
\begin{proposition}
\label{delta_proposition}
For all $\wt{x}\in\wt{\R^n}$ there exists $v\in\GS(\R^n)$ such that for all $u\in\Gtii(\R^n)$
\beq
\label{global_formula}
u(\wt{x})=\int_{\R^n}v(y)u(y)\, dy.
\eeq
\end{proposition}
\begin{proof}
We begin by observing that 
\begin{equation}
\label{tempered_v}
v:=(\varphi_\eps(x_\eps-\cdot))_\eps+\NS(\R^n)
\end{equation}
is a well-defined element of $\GS(\R^n)$. In fact since $(x_\eps)_\eps\in\R^n_M$ and $\varphi$ belongs to $\S(\R^n)$, for all $\alpha,\beta\in\N^n$ there exists $N\in\mathbb{N}$ depending on $\alpha$ such that the estimate
\begin{multline*}
|y^\alpha\partial^\beta_y(\varphi_\eps(x_\eps-y))|\le \eps^{-n-|\beta|}|y^\alpha\partial^\beta\varphi(\frac{x_\eps-y}{\eps})|\\
\le c_{1}\eps^{-n-|\beta|}\langle y\rangle^{|\alpha|}\langle\eps^{-1}({x_\eps-y})\rangle^{-|\alpha|}
\le c_2\eps^{-n-|\beta|}\langle x_\eps\rangle^{|\alpha|}\le c_3\eps^{-n-|\beta|-N}
\end{multline*}
is valid for $\eps$ varying in a small interval $(0,\eta]$. This shows that $(\varphi_\eps(x_\eps-\cdot))_\eps$ is a net belonging to $\ES(\R^n)$. A Taylor's formula argument together with Peetre's inequality shows that if $(x_\eps)_\eps$ and $(x'_\eps)_\eps$ are two representatives of $\tilde{x}$ then $(\varphi_\eps(x_\eps-y)-\varphi_\eps(x'_\eps-y))_\eps$ belongs to $\NS(\R^n)$.
%Let now $(x_\eps)_\eps$ and $(x'_\eps)_\eps$ be two representatives of $\tilde{x}$. We can write for $\eps$ small enough and arbitrary $q\in\mathbb{N}$,
%\begin{multline*}
%|y^\alpha\partial^\beta_y(\varphi_\eps(x_\eps-y)-\varphi_\eps(x'_\eps-y))|\le|y^\alpha\nabla(\partial^\beta\varphi_\eps)(x'_\eps-y+\theta(x_\eps-x'_\eps))||x_\eps-x'_\eps|\\
%\le\hskip-3pt c_1\eps^{-n-|\beta|-1}\langle y\rangle^{|\alpha|}\langle x'_\eps-y+\theta(x_\eps-x'_\eps)\rangle^{-|\alpha|}|x_\eps-x'_\eps|\le c_2\eps^{-n-|\beta|-1}\langle x'_\eps\rangle^{|\alpha|}|x_\eps-x'_\eps|\\
%\le c_3\eps^{-n-|\beta|-1-N|\alpha|+q},
%\end{multline*}
%where $N$ comes from the definition of $(x'_\eps)_\eps$. Thus $(\varphi_\eps(x_\eps-y)-\varphi_\eps(x'_\eps-y))_\eps\in\NS(\R^n)$.
Finally we prove that $(u_\eps(x_\eps)-\int_{\mathbb{R}^n}u_\eps(y)\varphi_\eps(x_\eps-y)\, dy)_\eps\in\Neg$.
By Taylor expansion of $u_\eps$ at $x_\eps$ we can write for $q\in\mathbb{N}$
\[
\int_{\mathbb{R}^n}u_\eps(y)\varphi_\eps(x_\eps-y)\, dy-u_\eps(x_\eps)=\sum_{|\alpha|=q+1}\frac{1}{\alpha !}\int_{\mathbb{R}^n}\partial^\alpha u_\eps(x_\eps-\eps\theta y)(-\eps y)^\alpha\varphi(y)\, dy.
\]
Hence when $\eps$ is close to $0$ we have that 
\begin{multline}
\label{final_estimate}
\biggl|u_\eps(x_\eps)-\int_{\mathbb{R}^n}u_\eps(y)\varphi_\eps(x_\eps-y)dy\biggr|\le c_1\eps^{-N+q+1}\,\cdot\\
\cdot\sum_{|\alpha|=q+1}\int_{\mathbb{R}^n}\langle x_\eps-\eps\theta y\rangle^N|y^\alpha\varphi(y)|dy\le c_2\eps^{-N+q+1-NM}
\end{multline}
Since $(u_\eps)_\eps\in\Et^{\infty,\infty}(\R^n)$ the exponent $N$ in \eqref{final_estimate} is independent of $q$ and $M$ characterizes the moderateness of $(x_\eps)_\eps$. This consideration completes the proof.
\end{proof}
\begin{remark}
\label{remark_gtii}
As for Theorem \ref{delta_theorem} we conclude with some comments on the result obtained above.
\begin{trivlist}
\item[(i)] Proposition \ref{delta_proposition} says that the restriction $\delta_{\wt{x}}\in\L(\Gt(\R^n),\wt{\C})$ to $\Gtii(\R^n)$ is an integral operator of the form $\int_{\R^n}v(y)\cdot\, dy$ with $v\in\GS(\R^n)$. Note that by the continuity of the map $\G^m_{N,\S}(\R^n)\to\wt{\C}:u\to\int_{\R^n}v(y)u(y)\, dy$ the integral operator $\int_{\R^n}v(y)\cdot dy$ belongs itself to $\L(\Gt(\R^n),\wt{\C})$.
\item[(ii)] The equality \eqref{global_formula} does not necessarily hold when $u$ is an element of $\Gt(\R^n)\setminus\Gtii(\R^n)$ and $v$ is given by \eqref{tempered_v}. In fact taking $\varphi\in\S(\R)$ with $\overline{\varphi}(0)-\Vert\varphi\Vert_2^2\neq 0$ as in $(ii)$ Remark \ref{remark_v} and $u=(\overline{\varphi}_\eps(-y))_\eps+\Nt(\R)\in\Gt(\R)\setminus\Gtii(\R)$, we conclude that $u(0)-\int_\R v(y)u(y)\, dy=[(1/\eps(\overline{\varphi}(0)-\Vert\varphi\Vert_2^2))_\eps]\neq 0$.
\end{trivlist}
\end{remark}
Combining reasonings analogous to the ones used in Theorem \ref{representation_theorem} with Proposition \ref{delta_proposition} and the point value theory in $\Gt(\R^n)$, we arrive at an integral representation for generalized functions in $\Gt^{\infty,\infty}(\R^n)$ by means of $v(x,y):=(\varphi_\eps(x-y))_\eps+\Nt(\R^{2n})$.
\begin{proposition}
\label{representation_prop}
There exists $v\in\Gt(\R^{2n})$ with a representative $(v_\eps)_\eps$ satisfying \eqref{condvrepr} such that 
\beq
\label{formula_representation_3}
u(x)=\int_{\R^n}v(x,y)u(y)\, dy\qquad \text{in}\ \Gt(\R^n)
\eeq
for all $u\in\Gtii(\R^n)$.
\end{proposition}

\section{Embedding theorems of Colombeau algebras into their duals}
This section is devoted to embedding the most common Colombeau algebras into their duals. We model our considerations on the classical results concerning $\D(\Om)$, $\S(\R^n)$, $\E(\Om)$ and the distributions spaces $\E'(\Om)$, $\S'(\R^n)$, $\D'(\Om)$, giving particular attention to some continuity and density issues. For the sake of simplicity and to the advantage of the reader we organize our reasoning in three subsections.
\subsection{Embeddings in the local context: $\G(\Om)\subseteq\L(\Gc(\Om),\wt{\C})$ and $\Gc(\Om)\subseteq\L(\G(\Om),\wt{\C})$}
In this subsection by the expression ``local context'' we mean to consider an open subset $\Om$ of $\R^n$ and the Colombeau algebras $\G(\Om)$, $\Ginf(\Om)$, $\Gc(\Om)$, $\Gc^\infty(\Om)$ constructed on it. We use the adjective ``local'' since the algebras which we will work with are characterized by representatives satisfying boundedness conditions on compact sets. We recall that $\G(\Om)$ and $\Ginf(\Om)$ are both Fr\'echet $\wt{\C}$-modules when topologized through the families of ultra-pseudo-seminorms $\{\mP_{K,j}\}_{K\Subset\Om, j\in\N}$ (\cite[Example 3.6]{Garetto:04b}) and $\{\mP_{\Ginf(K)}\}_{K\Subset\Om}$ (\cite[Example 3.12]{Garetto:04b}), respectively, while $\Gc(\Om)$ and $\Gc^\infty(\Om)$ are equipped with the strict inductive limit $\wt{\C}$-linear topology of the sequences $(\G_{K_n}(\Om), \{\mP_{\G_{K_n}(\Om),j}\}_{j\in\N})_{n\in\N}$ (\cite[Example 3.7]{Garetto:04b}) and $(\Ginf_{K_n}(\Om), \mP_{\Ginf_{K_n}(\Om)})_{n\in\N}$ (\cite[Example 3.13]{Garetto:04b}), respectively. Finally every topological dual like $\L(\Gc(\Om),\wt{\C})$, $\L(\G(\Om),\wt{\C})$ and so on is endowed with the $\wt{\C}$-linear topology $\beta_b$ of uniform convergence on bounded subsets. For this choice of topologies it is clear that the identity between representatives induces the following $\wt{\C}$-linear and continuous embeddings: $\Gc(\Om)\subseteq\G(\Om)$, $\Gc^\infty(\Om)\subseteq\Ginf(\Om)$, $\Ginf(\Om)\subseteq\G(\Om)$, $\Gc^\infty(\Om)\subseteq\Gc(\Om)$. Moreover, by the arguments presented before Theorem \ref{theorem_compact_support}, which concern $\Gc(\Om)$ and $\G(\Om)$ but can be easily adapted to $\Gc^\infty(\Om)$ and $\Ginf(\Om)$, we conclude that $\Gc(\Om)$ is dense in $\G(\Om)$ and $\Gc^\infty(\Om)$ is dense in $\Ginf(\Om)$.
\begin{theorem}
\label{theorem_local_embedding}
Integration defines a $\wt{\C}$-linear continuous embedding $u\to\big(v\to\int_\Om u(y)v(y)\, dy\big)$ of
\begin{itemize}
\item[(i)] $\G(\Om)$ into $\L(\Gc(\Om),\wt{\C})$,
\item[(ii)] $\Gc(\Om)$ into $\L(\G(\Om),\wt{\C})$.
\end{itemize}
Restriction defines a $\wt{\C}$-linear continuous embedding $T\to T_{\vert_{\cdot}}$ of 
\begin{itemize}
\item[(iii)] $\L(\G(\Om),\wt{\C})$ into $\L(\Gc(\Om),\wt{\C})$,
\item[(iv)] $\L(\Ginf(\Om),\wt{\C})$ into $\L(\Gc^\infty(\Om),\wt{\C})$.
\end{itemize}
\end{theorem}
\begin{proof}
We begin by observing that the injectivity of the integral operator $u\to\big(v\to\int_\Om u(y)v(y)\, dy\big)$ in $(i)$ and $(ii)$ follows from Proposition 2.11 in \cite{GGO:03}. By Theorem 1.26 in \cite{Garetto:04b} we know that if $B$ is a bounded subset of $\Gc(\Om)$ then it is contained in some $\G_K(\Om)$ and bounded there. This means that taking $K'\Subset\Om$ with $K\subseteq {\rm{int}}(K')$ we can write for all $u\in\G(\Om)$ and $v\in B$
\beq
\label{estimate_cont_embed}
\biggl\vert\int_\Om u(y)v(y)\, dy\biggr\vert_\esp=\biggl\vert\int_{K'} u(y)v(y)\, dy\biggr\vert_\esp \le \mP_{\G_K(\Om),0}(v)\, \mP_{K',0}(u)
\eeq
which leads to the estimate
\[
\sup_{v\in B}\biggl\vert\int_\Om u(y)v(y)\, dy\biggr\vert_\esp \le \big(\sup_{v\in B}\mP_{\G_K(\Om),0}(v)\big)\, \mP_{K',0}(u)\le C \mP_{K',0}(u)
\]
characterizing the continuity in $(i)$. Concerning $(ii)$ it suffices to show that for every $K\Subset\Om$ the map $\G_K(\Om)\to\L(\G(\Om),\wt{\C}):u\to\big(v\to\int_\Om u(y)v(y)\, dy\big)$ is continuous. When $B$ is a bounded subset of $\G(\Om)$ we immediately have from \eqref{estimate_cont_embed} that 
\[
\sup_{v\in B}\biggl\vert\int_\Om u(y)v(y)\, dy\biggr\vert_\esp \le \big(\sup_{v\in B}\mP_{K',0}(v)\big)\, \mP_{\G_K(\Om),0}(u)\le C \mP_{\G_K(\Om),0}(u).
\]
Finally the restriction maps in $(iii)$ and $(iv)$ are injective because of the density of $\Gc(\Om)$ in $\G(\Om)$ and of $\Gc^\infty(\Om)$ in $\Ginf(\Om)$. Moreover, by continuity arguments a bounded subset of $\Gc(\Om)$ is also bounded in $\G(\Om)$ and the same holds for the pairing $(\Gc^\infty(\Om),\Ginf(\Om))$. Therefore, the maps $\L(\G(\Om),\wt{\C})\to\L(\Gc(\Om),\wt{\C}): T\to T_{\vert_{\Gc(\Om)}}$ and $\L(\Ginf(\Om),\wt{\C})\to\L(\Gc^\infty(\Om),\wt{\C}): T\to T_{\vert_{\Gc^\infty(\Om)}}$ are continuous.
\end{proof}
Combining the results of Theorem \ref{theorem_local_embedding} we obtain the following chains of continuous inclusions
\beq
\label{chain_Ginf}
\Ginf(\Om)\subseteq\G(\Om)\subseteq\L(\Gc(\Om),\wt{\C}),
\eeq
\beq
\label{chain_Ginfc}
\Gc^\infty(\Om)\subseteq\Gc(\Om)\subseteq\L(\G(\Om),\wt{\C}),
\eeq
which play a relevant role in the regularity theory of differential and pseu\-do\-dif\-fe\-ren\-tial operators with generalized symbols (c.f. \cite{GGO:03}). In the following remarks we collect some further properties of the spaces involved in Theorem \ref{theorem_local_embedding}.
\begin{remark}
\label{remark_density}
\bf{Density of $\L(\G(\Om),\wt{\C})$ in $\L(\Gc(\Om),\wt{\C})$ and of $\L(\Ginf(\Om),\wt{\C})$ in $\L(\Gc^\infty(\Om),\wt{\C})$}\rm

$\L(\G(\Om),\wt{\C})$ is a dense subset of $\L(\Gc(\Om),\wt{\C})$ with respect to the topology of uniform convergence on bounded sets. In fact, choosing a partition of unity $(\chi_j)_{j\in\N}$ subordinate to a locally finite covering $(\Om_j)_j$ of $\Om$, for every $T\in\L(\Gc(\Om),\wt{\C})$ the map $T(\chi_j\,\cdot):\G(\Om)\to\wt{\C}:u\to T(\chi_ju)$ belongs to $\L(\G(\Om),\wt{\C})$. Further if $u\in\Gc(\Om)$ then $\sum_{j=0}^\infty \chi_ju$ is convergent to $u$ in $\Gc(\Om)$ since only a finite number of generalized functions $\chi_j u$ is different from $0$ when the support of $u$ is contained in some $K\Subset\Om$. Thus, by the properties of bounded subsets in $\Gc(\Om)$, the sequence $T_k:=\sum_{j=0}^k T(\chi_j\, \cdot)$ in $\L(\G(\Om),\wt{\C})$ tends to $T$ according to the topology $\beta_b(\L(\Gc(\Om),\wt{\C}),\Gc(\Om))$. Since the same proof remains valid for $\Gc^\infty(\Om)$ and $\Ginf(\Om)$ in place of $\Gc(\Om)$ and $\G(\Om)$ respectively, we have that $\L(\Ginf(\Om),\wt{\C})$ is a dense subset of $\L(\Gc^\infty(\Om),\wt{\C})$ endowed with $\beta_b(\L(\Gc^\infty(\Om),\wt{\C}),\Gc^\infty(\Om))$.
\end{remark}
\begin{remark}
\label{remark_local_embedding}
\bf{Failure of injectivity and density}\rm
\begin{trivlist}
\item[(i)] One may think of restricting $T\in\L(\G(\Om),\wt{\C})$ to the subalgebra of regular generalized functions $\Ginf(\Om)$. Obviously $T_{\vert_{\Ginf(\Om)}}$ belongs to $\L(\Ginf(\Om),\wt{\C})$ but the restriction map
\[
\L(\G(\Om),\wt{\C})\to\L(\Ginf(\Om),\wt{\C}): T\to T_{\vert_{\Ginf(\Om)}}
\]
is not injective. Theorem \ref{delta_theorem} allows to construct a map $T\neq 0$ in $\L(\G(\Om),\wt{\C})$ which is identically $0$ on $\Ginf(\Om)$. Take a mollifier $\varphi\in\S(\R^n)$ such that $\int_{\R^n}y_i\,\varphi^2(y)\, dy \neq 0$ for some $i\in\{1,...,n\}$. Let $x_0$ be an arbitrary point of $\Om$ and $\psi$ a cut-off function in $\Cinfc(\Om)$ equal to $1$ in a neighborhood of $x_0$. The generalized function $(\psi(y)(x_{0,i}-y_i)\varphi_\eps(x_0-y))_\eps +\Neg(\Om)$ belongs to $\Gc(\Om)$ and thus
\beq
\label{definition_T}
T:\G(\Om)\to\wt{\C}:u\to \biggl(\int_\Om u_\eps(y)(x_{0,i}-y_i)\psi(y)\varphi_\eps(x_0-y)\, dy\biggr)_\eps +\Neg
\eeq
is an element of $\L(\G(\Om),\wt{\C})$. The restriction of $T$ to $\Ginf(\Om)$ is identically $0$ since when $u$ belongs to $\Ginf(\Om)$ the integral in \eqref{definition_T} is the point value of $u(y)(x_{0,i}-y_i)_{\vert_{y\in\Om}}$ at $y=x_0$. Consider now the generalized function $u_0:=[(\varphi_\eps(x_0-\cdot)_{\vert_\Om})_\eps]$ in $\G(\Om)\setminus\Ginf(\Om)$. $Tu_0$ has the net
\begin{multline}
\label{rapresent_Tu}
\int_\Om \psi(y)(x_{0,i}-y_i)(\varphi_\eps)^2(x_0-y)\, dy = \int_{\R^n}(x_{0,i}-y_i)(\varphi_\eps)^2(x_0-y)\, dy\\ + \int_{\R^n}(\psi(y)-1)(x_{0,i}-y_i)(\varphi_\eps)^2(x_0-y)\, dy
\end{multline}
as a representative. Computations analogous to \eqref{rap_decreasing_varphi} show that the second summand in the right-hand side of \eqref{rapresent_Tu} is negligible. Hence $Tu_0=[(\eps^{1-n}\int_{\R^n}y_i\varphi^2(y)\, dy)_\eps]\neq 0$. This same example shows that the restriction map
\[
\L(\G(\Om),\wt{\C})\to\L(\Gc^\infty(\Om),\wt{\C}): T\to T_{\vert_{\Gc^\infty(\Om)}}
\]
is not injective. Finally an easy adaptation of the previous arguments to the pairing $(\Gcinf(\Om),\Gc(\Om))$ combined with $(iii)$ Remark \ref{remark_v} proves the non-injectivity of $\L(\Gc(\Om),\wt{\C})\to\L(\Gc^\infty(\Om),\wt{\C}): T\to T_{\vert_{\Gc^\infty(\Om)}}$. As a consequence of these results we obtain that $\Ginf(\Om)$ and $\Gc^\infty(\Om)$ are not dense in $\G(\Om)$ and $\Gc(\Om)$ respectively.
%Moreover defining $T\in\L(\Gc(\Om),\wt{\C})$ as the map which assigns the integral $[(\int_\Om (x_{0,i}-y_i)u_\eps(y)\varphi_\eps(x_0-y)\, dy)_\eps]$ to $u\in\Gc(\Om)$ ($\supp(u_\eps)\subseteq K\Subset\Om$ for all $\eps\in(0,1]$), we have by $(iv)$Remark \ref{remark_v} that  $Tu=0$ for all $u\in\Gc^\infty(\Om)$ while for $u_0=(\psi(y)\varphi_\eps(x_0-y))_\eps +\Neg(\Om)\in\Gc(\Om)\setminus\Gc^\infty(\Om)$ we have that $Tu_0\neq 0$. In conclusion this proves the non-injectivity of the restriction
%\[
%\L(\Gc(\Om),\wt{\C})\to\L(\Gc^\infty(\Om),\wt{\C}): T\to T_{\vert_{\Gc^\infty(\Om)}}.
%\]
%\item[(ii)] We already observed that $\Gc(\Om)$ is dense in $\G(\Om)$ and that $\Gc^\infty(\Om)$ is dense in $\Ginf(\Om)$. As a consequence of the previous discussion we obtain that $\Ginf(\Om)$ and $\Gc^\infty(\Om)$ are not dense in $\G(\Om)$ and $\Gc(\Om)$ respectively.
\item[(ii)] Every $u\in\G(\Om)$ defines an element of $\L(\Gc^\infty(\Om),\wt{\C})$ via $\int_\Om u(y)\cdot\, dy$. However, $\G(\Om)$ is not contained in $\L(\Gc^\infty(\Om),\wt{\C})$ since the map 
\[
\G(\Om)\to\L(\Gc^\infty(\Om),\wt{\C}):u\to \biggl(v\to\int_\Om u(y)v(y)\, dy\biggr)
\]
is not injective. In fact as shown in the first assertion of this remark, for $u=[((x_{0,i}-y_i)_{\vert_{y\in\Om}}\varphi_\eps(x_0-y)_{\vert_{y\in\Om}})_\eps]\in\G(\Om)$ the integral $\int_\Om u(y)v(y)\, dy$ is $0$ for all $v\in\Gc^\infty(\Om)$ but $u$ is not $0$ in $\G(\Om)$. Applying the same kind of reasoning to $u=[(\psi(y)(x_{0,i}-y_i)\varphi_\eps(x_0-y))_\eps]\in\Gc(\Om)$ we obtain that the map
\[
\Gc(\Om)\to\L(\Ginf(\Om),\wt{\C}):u\to \biggl(v\to\int_\Om u(y)v(y)\, dy\biggr)
\]
is not injective.
\end{trivlist}
\end{remark}
%We conclude this subsection by discussing the following topological issue.
%By \cite[Propositions 2.9, 2.14, 2.15]{Garetto:04b} the Colombeau algebras $\Ginf(\Om)$, $\Gc^\infty(\Om)$, $\GSinf(\R^n)$ are all bornological and barrelled locally convex topological $\wt{\C}$-modules. Hence their topological duals are complete for the corresponding strong topology and for the topology $\beta_b$ of uniform convergence on compact sets. In addition endowed with the weak topology these duals are quasi-complete. Theorem 1.25 in \cite{Garetto:04b} applied to $\Gc^\infty(\Om)$ says that a subset $B$ of $\Gc^\infty(\Om)$ is bounded if and only if it is contained in some $\G^\infty_K(\Om)$ and bounded there. Finally we have that Corollary 2.18 in \cite{Garetto:04b} can be stated for $\G(\Om)$, $\Gc(\Om)$, $\GS(\R^n)$ as well as for their regular subalgebras $\Ginf(\Om)$, $\Gc^\infty(\Om)$, $\GSinf(\R^n)$.
\subsection{Embedding in the global context: $\GS(\R^n)\subseteq\Gt(\R^n)\subseteq\L(\GS(\R^n),\wt{\C})$}
The expression ``global context'' emphasizes the fact that the algebras considered in this subsection are given via global estimates on $\R^n$ at the level of representatives. 

Our aim is now to obtain chains of inclusions similar to \eqref{chain_Ginf} and \eqref{chain_Ginfc} which will involve $\GSinf(\R^n)$, $\GS(\R^n)$ and $\Gt(\R^n)$. In this procedure some auxiliary algebras of generalized functions, whose properties turn out to be useful in many situations, are introduced and some preliminary results are discussed. More precisely we refer to the algebras $\Gpq(\R^n)$ from \cite{BO:92} and $\G^\infty_{p,q}(\R^n)$, $1\le p,q\le \infty$, defined as the quotients $\EMp(\R^n)/\Npq(\R^n)$ and $\EMp^\infty(\R^n)/\Npq(\R^n)$ respectively, where
\beq
\label{def_Gpq}
\begin{split}
\E_p[\R^n]& := W^{\infty,p}(\R^n)^{(0,1]},\\
\EMp(\R^n)& :=\M_{W^{\infty,p}(\R^n)},\\
\EMp^\infty(\R^n) & := \M^\infty_{W^{\infty,p}(\R^n)},\\
\Npq(\R^n) & :=\{(u_\eps)_\eps\in\EMp(\R^n)\cap\E_q[\R^n]:\\
&\qquad\qquad\qquad\forall\alpha\in\N^n\ \forall m\in\N\quad\Vert\partial^\alpha u_\eps\Vert_q= O(\eps^m)\ \text{as}\ \eps\to 0\}.
\end{split}
\eeq
The space $W^{-\infty,p}(\R^n)$ is embedded into $\Gpq(\R^n)$ via convolution and for $1\le p_1\le p_2\le\infty$, $1\le q\le\infty$ we have that $\G_{p_1,q}(\R^n)\subseteq\G_{p_2,q}(\R^n)$. Note that $\Gpp(\R^n)=\G_{W^{\infty,p}(\R^n)}$ (\cite[Definition 3.1]{Garetto:04b}) and $\Gpp^\infty(\R^n)=\G^\infty_{W^{\infty,p}(\R^n)}$ (\cite[Example 3.10]{Garetto:04b}). Therefore, equipped with the corresponding sharp topologies, $\Gpp(\R^n)$ is a Fr\'echet $\wt{\C}$-module and $\Gpp^\infty(\R^n)$ is a complete ultra-pseudo-normed $\wt{\C}$-module.

We begin by investigating the structure of the ideals $\NS(\R^n)$ and $\Npp(\R^n)$ pointing out the algebraic relationships between $\GS(\R^n)$ and $\Gpq(\R^n)$. In the sequel $\ES(\R^n)$ denotes the space $\M_{\S(\R^n)}$ of moderate nets and we slightly simplify the general expressions of Section 3 in \cite{Garetto:04b} concerning $\G_E$ by calling the elements of $\NS(\R^n)$ and $\Npp(\R^n)$ $\S$-negligible and $p$-negligible, respectively.
\begin{proposition}
\label{prop_ideal_NS_Npp}
\leavevmode
\begin{itemize}
\item[(i)] $(u_\eps)_\eps\in\ES(\R^n)$ is $\S$-negligible if and only if the following condition is satisfied:
\beq
\label{charns}
\forall q\in\mathbb{N}\qquad\qquad \sup_{x\in\mathbb{R}^n}|u_\eps(x)|=O(\eps^q)\qquad as\ \eps\to 0.
\eeq
\item[(ii)] $(u_\eps)_\eps\in\EMp(\R^n)$ is $p$-negligible if and only if 
\beq
\label{charnp}
\forall m\in\N\qquad\qquad \Vert u_\eps\Vert_p=O(\eps^m)\qquad as\ \eps\to 0.
\eeq
\end{itemize}
\end{proposition}
These statements allow to add $\NS(\R^n)$ and $\Npp(\R^n)$ to the list of ideals of a Colombeau algebra which, once moderateness is known, can be characterized by estimates on the $0$-th derivative. This list already contains $\Neg(\Om)$ (Theorem 1.2.3 \cite{GKOS:01}) and $\Nt(\R^n)$ (Theorem 1.2.25 \cite{GKOS:01}).
\begin{proof}
$(i)$ It is clear that \eqref{charns} holds if $(u_\eps)_\eps$ is $\S$-negligible. Conversely assume $(u_\eps)_\eps\in\ES(\R^n)$ and that \eqref{charns} is valid. We shall show that 
\begin{multline}
\label{1step}
\forall q\in\mathbb{N}\ \forall N\in\mathbb{N}\ \exists\eta\in(0,1]\ \forall x\in\mathbb{R}^n\  \forall\eps\in(0,\eta]\\
(1+|x|)^{N}|u_\eps(x)|\le \eps^q
\end{multline}
and using \eqref{1step} that
\begin{multline}
\label{2step}
\forall q\in\mathbb{N}\ \forall N\in\mathbb{N}\ \forall i=1,...,n,\ \exists\eta\in(0,1]\ \forall x\in\mathbb{R}^n\ \forall\eps\in(0,\eta]\\
\quad(1+|x|)^N|\partial_iu_\eps(x)|\le \eps^q.
\end{multline}
The same reasoning applied to the derivatives of order one allows us to extend the result to derivatives of any order, obtaining the estimates of an $\S$-negligible net. We begin then by proving \eqref{1step}. Since $(u_\eps)_\eps\in\ES(\R^n)$ we know that 
\begin{equation}
\label{mods1}
\forall N\in\mathbb{N}\ \exists M\in\mathbb{N}\qquad\quad \sup_{x\in\mathbb{R}^n}(1+|x|)^N|u_\eps(x)|=O(\eps^{-M})\quad \text{as}\ \eps\to 0.
\end{equation}
By \eqref{charns}, for any $q\in\mathbb{N}$ there exists $\eta\in(0,1]$ small enough, such that for all $x\in\mathbb{R}^n$ and $\eps\in(0,\eta]$ we have
\beq
\label{est1ns}
(1+|x|)^N|u_\eps(x)|^2=(1+|x|)^N|u_\eps(x)||u_\eps(x)|\le \eps^{q-M},
\eeq
where $M$ depends on $N$ as in \eqref{mods1}. Since the choice of $N$ and $q$ in \eqref{est1ns} is arbitrary, this estimate leads to \eqref{1step}.

We turn to the proof of \eqref{2step}. From the assumption of $\S$-moderateness and the estimate \eqref{1step} we have for $i=1,...,n$ that
\beq
\label{est2ns}
\forall N\in\mathbb{N}\ \exists M\in\mathbb{N}\qquad \sup_{x\in\mathbb{R}^n}(1+|x|)^N|\partial^2_iu_\eps(x)|=O(\eps^{-M})\quad \text{as}\ \eps\to 0
\eeq
and
\beq
\label{est3ns}
\forall q\in\mathbb{N}\qquad\qquad\sup_{x\in\mathbb{R}^n}(1+|x|)^N|u_\eps(x)|=O(\eps^{2q+M})\quad \text{as}\ \eps\to 0.
\eeq
Fix now $N\in\mathbb{N}$. By Taylor's formula we can write
\[
u_\eps(x+\eps^{M+q}e_i)=u_\eps(x)+\partial_iu_\eps(x)\eps^{M+q}+\frac{1}{2}\partial^2_iu_\eps(x+\eps^{M+q}\theta e_i)\eps^{2M+2q},
\]
where $e_i$ is the i-th vector of the canonical basis of $\mathbb{R}^n$ and $\theta\in[0,1]$. Consequently,
\begin{multline*}
(1+|x|)^N\partial_iu_\eps(x)=(1+|x|)^N(u_\eps(x+\eps^{M+q}e_i)-u_\eps(x))\eps^{-M-q}+\\
-\frac{1}{2}(1+|x|)^N\partial^2_iu_\eps(x+\eps^{M+q}\theta e_i)\eps^{M+q}.
\end{multline*}
Since $(1+|x|)^N=(1+|x+\eps^{M+q}e_i-\eps^{M+q}e_i|)^N\le(1+|x+\eps^{M+q}e_i|+\eps^{M+q})^N$, we conclude that 
\beq
\label{est4ns}
\begin{split}
\sup_{x\in\mathbb{R}^n}(1+|x|)^N &|\partial_iu_\eps(x)|\le  2^N\eps^{-M-q}\sup_{x\in\mathbb{R}^n}(1+|x+\eps^{M+q}e_i|)^N|u_\eps(x+\eps^{M+q}e_i)|\\
&+\eps^{-M-q}\sup_{x\in\mathbb{R}^n}(1+|x|)^N|u_\eps(x)|\\
&+2^{N-1}\eps^{M+q}\sup_{x\in\mathbb{R}^n}(1+|x+\theta\eps^{M+q}e_i|)^N|\partial^2_iu_\eps(x+\eps^{M+q}\theta e_i)|.
\end{split}
\eeq
By \eqref{est3ns} and \eqref{est2ns}, the first two suprema over $\mathbb{R}^n$ in the right-hand side of \eqref{est4ns} are $O(\eps^{2q+M})$ as $\eps\to 0$ and the third is $O(\eps^{-M})$ as $\eps\to 0$. These considerations imply $\sup_{x\in\mathbb{R}^n}(1+|x|)^N|\partial_iu_\eps(x)|=O(\eps^q)$ and complete the proof of $(i)$.

$(ii)$ This assertion is obtained from the proof of Theorem 1.2.3 in \cite{GKOS:01} if we take global estimates on $\R^n$ and we substitute the $L^\infty$-norm with the $L^p$-norm.
\end{proof}
Proposition \ref{prop_ideal_NS_Npp} makes it easy to prove the following injectivity result.
\begin{proposition}
\label{prop_GS_Gpq}
The map
\beq
\label{map_GS_Gpq}
\GS(\R^n)\to\Gpq(\R^n):(u_\eps)_\eps+\NS(\R^n)\to(u_\eps)_\eps+\Npq(\R^n)
\eeq
is well-defined, injective and maps $\GSinf(\R^n)$ into $\Gpq^\infty(\R^n)$.
\end{proposition}
\begin{proof} 
The well-definedness of \eqref{map_GS_Gpq} and the additional property $\GSinf(\R^n)\to\Gpq^\infty(\R^n)$ come from the fact that $\ES(\R^n)\subseteq\EMp(\R^n)$, $\ES^\infty(\R^n)\subseteq\EMp^\infty(\R^n)$ and $\NS(\R^n)\subseteq\Npq(\R^n)$. It remains to show that $(u_\eps)_\eps\in\ES(\R^n)\cap\Npq(\R^n)$ implies $(u_\eps)_\eps\in\NS(\R^n)$. Since it suffices to estimate the $L^\infty$-norm of $(u_\eps)_\eps$, we conclude by Sobolev's embedding theorem that
\[
\Vert u_\eps\Vert_\infty\le c\max_{\substack{|\beta|\le k,\\ kq>n}}\Vert\partial^\beta u_\eps\Vert_q=O(\eps^m)\ \ \text{as}\ \eps\to 0,
\]
for arbitrary $m\in\mathbb{N}$ .
\end{proof}
In the case of $p=q$, Proposition \ref{prop_GS_Gpq} and the continuity of the embedding $\S(\R^n)\subseteq W^{\infty,p}(\R^n)$ applied at the level of  representatives, yields continuity with respect to the sharp topologies of the inclusions $\GS(\R^n)\subseteq\Gpp(\R^n)$ and $\GSinf(\R^n)\subseteq\Gpp^\infty(\R^n)$. 

In the next proposition we establish a topological link between $\Gc(\Om)$ and $\Gpp(\R^n)$. The following embeddings are defined using representatives $(u_\eps)_\eps$ of $u\in\Gc(\Om)$ with $\supp\, u_\eps$ contained in some compact set of $\Om$ uniformly with respect to $\eps$. 
\begin{proposition}
\label{prop_chain_Gc_Gpp}
Let $\Om$ be an open subset of $\R^n$. The following chains of $\wt{\C}$-linear continuous embeddings between locally convex topological $\wt{\C}$-modules hold:
\beq
\label{chain_Gc_Gpp}
\Gc(\Om)\subseteq\Gc(\R^n)\subseteq\GS(\R^n)\subseteq\Gpp(\R^n),
\eeq
\beq
\label{chain_Gcinf_Gppinf}
\Gcinf(\Om)\subseteq\Gcinf(\R^n)\subseteq\GSinf(\R^n)\subseteq\G_{p,p}^\infty(\R^n)
\eeq
\end{proposition}
%\begin{tabular}{ccccccc}
%$\Gc(\Om)\hskip-10pt$ & $\hskip-10pt\longrightarrow\hskip-10pt$ & $\hskip-10pt\Gc(\R^n)\hskip-10pt$ & $\hskip-10pt\longrightarrow\hskip-10pt$ & $\hskip-10pt\GS(\R^n)\hskip-10pt$ & $\hskip-10pt\longrightarrow\hskip-10pt$ & $\hskip-10pt\Gpp(\R^n)$\\[0.3cm]
%$(u_\eps)_\eps\hskip-3pt+\hskip-3pt\Neg(\Om)\hskip-10pt$ & ${\, }$ & $\hskip-10pt(u_\eps)_\eps\hskip-3pt+\hskip-3pt\Neg(\R^n)\hskip-10pt$ & ${\, }$ & $\hskip-10pt(u_\eps)_\eps\hskip-3pt+\hskip-3pt\NS(\R^n)\hskip-10pt$ & ${\, }$ & $\hskip-10pt(u_\eps)_\eps\hskip-3pt+\hskip-3pt\Npp(\R^n)$
%\end{tabular}
%\eeq
%and
%\beq
%\label{chain_Gcinf_Gppinf}
%\begin{tabular}{ccccccc}
%$\Gc^\infty(\Om)\hskip-10pt$ & $\hskip-10pt\longrightarrow\hskip-10pt$ & $\hskip-10pt\Gc^\infty(\R^n)\hskip-10pt$ & $\hskip-10pt\longrightarrow\hskip-10pt$ & $\hskip-10pt\GSinf(\R^n)\hskip-10pt$ & $\hskip-10pt\longrightarrow\hskip-10pt$ & $\hskip-10pt\Gpp^\infty(\R^n)$\\[0.3cm]
%$(u_\eps)_\eps\hskip-3pt+\hskip-3pt\Neg(\Om)\hskip-10pt$ & ${\, }$ & $\hskip-10pt(u_\eps)_\eps\hskip-3pt+\hskip-3pt\Neg(\R^n)\hskip-10pt$ & ${\, }$ & $\hskip-10pt(u_\eps)_\eps\hskip-3pt+\hskip-3pt\NS(\R^n)\hskip-10pt$ & ${\, }$ & $\hskip-10pt(u_\eps)_\eps\hskip-3pt+\hskip-3pt\Npp(\R^n)$
%\end{tabular}
%\eeq
\begin{remark}
\label{remark_consequences}
The apparently obvious claim of Proposition \ref{prop_chain_Gc_Gpp} has some interesting consequences at the level of topological duals. First of all it provides an alternative to Proposition 2.11 in \cite{GGO:03} for proving the inclusion of $\Gc(\Om)$ in $\L(\G(\Om),\wt{\C})$ and of $\G(\Om)$ in $\L(\Gc(\Om),\wt{\C})$. Indeed if we assume $u\in\Gc(\Om)$ and $\int_\Om u(y)v(y)\, dy=0$ for all $v\in\G(\Om)$, then taking $v=\overline{u}:=(\overline{u_\eps})_\eps+\Neg(\Om)$ we have that $\int_\Om |u|^2(y)\, dy=0$. This means that $u$ is $0$ in $\G_{2,2}(\R^n)$ and therefore by \eqref{chain_Gc_Gpp} it is $0$ in $\Gc(\Om)$. Analogously when $u$ belongs to $\G(\Om)$ and it is the $0$-element of $\L(\Gc(\Om),\wt{\C})$, its support is forced to be empty since $u\psi=0$ in $\G_{2,2}(\R^n)$ for all cut-off functions $\psi$. Clearly using \eqref{chain_Gcinf_Gppinf} we can prove that the integral operator $u\to (v\to\int_\Om u(y)v(y)\, dy)$ gives a $\wt{\C}$-linear embedding of $\Gc^\infty(\Om)$ into $\L(\Ginf(\Om),\wt{\C})$ and of $\Ginf(\Om)$ into $\L(\Gc^\infty(\Om),\wt{\C})$.
\end{remark}
We are ready now for dealing with $\GS(\R^n)$, $\Gt(\R^n)$ and the dual $\L(\GS(\R^n),\wt{\C})$ stating the ``global version'' of Theorem \ref{theorem_local_embedding}. We recall that $\GS(\R^n)$ is a Fr\'echet $\wt{\C}$-module if topologized through the family of ultra-pseudo-seminorms $\{\mP_j\}_{j\in\N}$ (\cite[Example 3.6]{Garetto:04b}) while ($\GSinf(\R^n),\mP^\infty_{\S(\R^n)})$ is a complete ultra-pseudo-normed $\wt{\C}$-module (\cite[Example 3.11]{Garetto:04b}). Note that for this choice of topologies $\GSinf(\R^n)$ is continuously embedded in $\GS(\R^n)$. The topology on $\Gt(\R^n)$ is defined via $\GtS(\R^n)=\cap_{m\in\R}\cup_{N\in\N}\G^m_{N,\S}(\R^n)$ (\cite[Example 3.9]{Garetto:04b}) where on each $\G^m_{N,\S}(\R^n)$ we take the ultra-pseudo-seminorm $\mP^m_N$. Finally we equip $\L(\GS(\R^n),\wt{\C})$ with the topology of uniform convergence on bounded subsets of $\GS(\R^n)$.
\begin{theorem}
\label{theorem_global_embedding}
Identity between representatives defines a $\wt{\C}$-linear continuous embedding of 
\begin{itemize}
\item[(i)] $\GS(\R^n)$ into $\Gt(\R^n)$.
\end{itemize}
Integration defines a $\wt{\C}$-linear continuous embedding $u\to\big(v\to\int_\Om u(y)v(y)\, dy\big)$ of
\begin{itemize}
\item[(ii)] $\Gt(\R^n)$ into $\L(\GS(\R^n),\wt{\C})$.
\end{itemize}
\end{theorem}
\begin{proof}
$(i)$ Since $\ES(\R^n)\subseteq\Et(\R^n)$ and $\NS(\R^n)\subseteq\Nt(\R^n)$, the map $\GS(\R^n)\to\Gt(\R^n):(u_\eps)_\eps+\NS(\R^n)\to (u_\eps)_\eps+\Nt(\R^n)$ is well-defined. In addition it is injective since $\ES(\R^n)\cap\Nt(\R^n)\subseteq\NS(\R^n)$. In fact if $(u_\eps)_\eps$ is $\tau$-negligible then 
\[
\exists N\in\mathbb{N}\ \forall q\in\mathbb{N}\qquad\qquad \sup_{x\in\mathbb{R}^n}(1+|x|)^{-N}|u_\eps(x)|=O(\eps^q)\quad \text{as}\ \eps\to 0
\]
and by the $\S$-moderateness there exists $M\in\mathbb{N}$ depending on $N$ such that  $\sup_{x\in\mathbb{R}^n}(1+|x|)^{N}|u_\eps(x)|=O(\eps^{-M})$. Hence for any $q\in\N$, $x\in\R^n$ and $\eps$ small enough 
$$|u_\eps(x)|^2=(1+|x|)^{-N}|u_\eps(x)|(1+|x|)^{N}|u_\eps(x)|\le \eps^{q-M}$$
which means that $(u_\eps)_\eps\in\NS(\R^n)$ by Proposition \ref{prop_ideal_NS_Npp}. Note that the continuity in $(i)$ is guaranteed because $\GS(\R^n)$ is continuously embedded in every $\G^m_{N,\S}(\R^n)$.

$(ii)$ Concerning the second assertion, we first consider the continuity issue and then we deal with injectivity. If $u\in\G^m_{N,\S}(\R^n)$ then for all $v\in\GS(\R^n)$ we can write at the level of representatives
\begin{multline*}
\biggl|\int_{\R^n}u_\eps(y)v_\eps(y)\, dy\biggr|\le \sup_{y\in\R^n,|\alpha|\le m}(1+|y|)^{-N}|\partial^\alpha u_\eps(y)|\ \int_{\R^n}(1+|y|)^N|v_\eps(y)|\, dy\\
\le C\sup_{y\in\R^n,|\alpha|\le m}(1+|y|)^{-N}|\partial^\alpha u_\eps(y)|\ \sup_{y\in\R^n}(1+|y|)^{n+1+N}|v_\eps(y)|.
\end{multline*}
This says that for all $u\in\G^m_{N,\S}(\R^n)$ the integral $\int_{\R^n}u(y)\cdot\, dy$ is a map in $\L(\GS(\R^n),\wt{\C})$ and when $B$ is a bounded subset of $\GS(\R^n)$ yields
\beq
\label{est_cont_Gt}
\sup_{v\in B}\biggl\vert\int_{\R^n}u(y)v(y)\, dy\biggr\vert_\esp \le (\sup_{v\in B}\mP_{N+n+1}(v))\, \mP^m_N(u).
\eeq 
By \eqref{est_cont_Gt} and the topology defined on $\Gt(\R^n)$ we conclude that $u\to (v\to\int_{\R^n}u(y)v(y)\, dy)$ is a continuous map between $\Gt(\R^n)$ and $\L(\GS(\R^n),\wt{\C})$.

Now we want to prove that if $\int_{\R^n}u(y)v(y)\, dy=0$ for all $v\in\GS(\R^n)$ then $u$ is $0$ in $\Gt(\R^n)$. Fix $\wt{x}\in\wt{\R^n}$, a representative $(x_\eps)_\eps$ of $\wt{x}$ with $|x_\eps|=O(\eps^{-N})$ for some $N\in\N$ and a mollifier $\varphi\in\S(\R^n)$. Choosing a representative $(u_\eps)_\eps$ of $u$ the net $(\overline{u}_\eps|\widehat{\varphi}(\eps^{N+1}\cdot)|^2)_\eps$ belongs to $\ES(\R^n)$ since the estimate
\[
\begin{split}
|x^\alpha\partial^\beta(\overline{u}_\eps(x)|\widehat{\varphi}(\eps^{N+1}x)|^2)|&\le\big|x^\alpha\sum_{\gamma\le\beta}\binom{\beta}{\gamma}\partial^\gamma\overline{u}_\eps(x)\partial^{\beta-\gamma}|\widehat{\varphi}(\eps^{N+1}x)|^2\big|\\
&\le c(\eps^{N+1})^{-M}(1+|x|)^{|\alpha|-M}\sum_{\gamma\le\beta}|\partial^\gamma\overline{u}_\eps(x)|,
\end{split}
\]
is valid for all $M\in\N$. Thus $v:=(\overline{u}_\eps|\widehat{\varphi}(\eps^{N+1}\cdot)|^2)_\eps +\NS(\R^n)$ is a well-defined element of $\GS(\R^n)$ and by assumption $\int_{\R^n}u(y)v(y)\, dy=0$. This means that $(\Vert u_\eps\widehat{\varphi_{\eps^{N+1}}}\Vert_2^2)_\eps\in\Neg$ or in other words that $( u_\eps\widehat{\varphi_{\eps^{N+1}}})_\eps\in\ES(\R^n)\cap\Neg_{2,2}(\R^n)$ by $(ii)$, Proposition \ref{prop_ideal_NS_Npp}. As a consequence of Proposition \ref{prop_GS_Gpq} we have that $( u_\eps\widehat{\varphi_{\eps^{N+1}}})_\eps\in\NS(\R^n)$. Recalling that $\widehat{\varphi}(0)=1$, there exist some constants $c>0$ and $\eta\in(0,1]$ such that $|\widehat{\varphi}(x)|\ge 1/2$ for $|x|\le c$ and $|\eps^{N+1}x_\eps|\le c$ for all $\eps\in(0,\eta]$. Hence by the $\S$-negligibility-estimate characterizing $( u_\eps\widehat{\varphi_{\eps^{N+1}}})_\eps$ we conclude that for all $q\in\N$ there exists $\eta'\in(0,1]$ such that 
\[
|u_\eps(x_\eps)|=|u_\eps(x_\eps)\widehat{\varphi}(\eps^{N+1}x_\eps)|\frac{1}{|\widehat{\varphi}(\eps^{N+1}x_\eps)|}\le 2|u_\eps(x_\eps)\widehat{\varphi}(\eps^{N+1}x_\eps)|\le 2\eps^q.
\]
on the interval $(0,\eta']$. This shows that $u(\tilde{x})=0$ in $\wt{\C}$ and since $\wt{x}$ was arbitrary we have $u=0$ in $\Gt(\R^n)$.
\end{proof}
Note that by the continuity of the inclusion $\Gt(\R^n)\subseteq\L(\GS(\R^n),\wt{\C})$ and the separatedness of $\L(\GS(\R^n),\wt{\C})$ it follows that $\Gt(\R^n)$ is a separated locally convex topological $\wt{\C}$-module.

Summarizing Theorem \ref{theorem_global_embedding} implies the chain of continuous embeddings
\beq
\label{chain_GSinf_LGS}
\GSinf(\R^n)\subseteq \GS(\R^n)\subseteq \Gt(\R^n)\subseteq \L(\GS(\R^n),\wt{\C}).
\eeq
This systematic framework of topological $\wt{\C}$-modules is the core of the regularity theory for pseudodifferential operators with global generalized symbols \cite{Garetto:04th}.
\begin{remark}
\label{delta_Gs_remark}
It is clear from \eqref{chain_GSinf_LGS} that the algebra $\GS(\R^n)$ inherits the point value theory from $\Gt(\R^n)$. Hence it is meaningful for all $\wt{x}\in\wt{\R^n}$ to define the generalized delta functional $\delta_{\wt{x}}$ as an element of $\L(\GS(\R^n),\wt{\C})$. Proposition \ref{delta_proposition} can be applied to $u\in\GSinf(\R^n)$ and states that the restriction of $\delta_{\wt{x}}$ to $\GSinf(\R^n)$ is an integral operator of the form $\int_{\R^n}v(y)\cdot\, dy$, where $v\in\GS(\R^n)$. Note that by $(ii)$ Remark \ref{remark_gtii} the generalized function $v\in\GS(\R^n)$ defined in \eqref{tempered_v} does not necessarily provide the equality \eqref{global_formula} for $u\in\GS(\R^n)\setminus\GSinf(\R^n)$. Finally Proposition \ref{representation_prop} gives an integral representation for every $u\in\GSinf(\R^n)$. A straightforward application of Proposition \ref{delta_proposition} allows us to adapt $(i)$ Remark \ref{remark_local_embedding} to the context of $\GSinf(\R^n)$ and $\GS(\R^n)$, proving that the following restriction
\[
\L(\GS(\R^n),\wt{\C})\to\L(\GSinf(\R^n),\wt{\C}):T\to T_{\vert_{\GSinf(\R^n)}}
\]
is continuous but not injective. Consequently $\GSinf(\R^n)$ is not dense in $\GS(\R^n)$ and the algebra $\GS(\R^n)$ is not contained in $\L(\GSinf(\R^n),\wt{\C})$ since the $\wt{\C}$-linear map $\GS(\R^n)\to\L(\GSinf(\R^n),\wt{\C}):u\to(v\to\int_{\R^n}u(y)v(y)\, dy)$ is not an injection.
\end{remark}
\subsection{Distribution theory and duality theory for Colombeau algebras: some comparisons and a regularization of $\delta_{\wt{x}}$}
We conclude the paper by discussing some issues which relate distribution theory to duality theory for Colombeau algebras. 

First we concentrate on embedding procedures. Example \ref{example_distributions} provides a simple way of embedding distributions into the topological dual of a Colombeau algebra, which essentially consists in using the distribution as a representative of a $\wt{\C}$-linear and continuous map. Now that the chains of inclusions $\E'(\Om)\subseteq\Gc(\Om)\subseteq\L(\G(\Om),\wt{\C})$, $\D'(\Om)\subseteq\G(\Om)\subseteq\L(\Gc(\Om),\wt{\C})$ and $\S'(\R^n)\subseteq\Gt(\R^n)\subseteq\L(\GS(\R^n),\wt{\C})$ provide us with an alternative way of interpreting distributions as elements of the dual of a suitable Colombeau algebra, we would like to compare the two types of embedding. We fix and recall some notations. We denote the straightforward embedding of $\D'(\Om)$ into $\L(\Gc(\Om),\wt{\C})$ introduced in Example \ref{example_distributions} by $\iota_d$ and the embedding of $\D'(\Om)$ into the Colombeau algebra $\G(\Om)$ by $\iota$. Finally we denote the composition of $\iota$ with the embedding of $\G(\Om)$ into the dual $\L(\Gc(\Om),\wt{\C})$ (Theorem \ref{theorem_local_embedding}) by $\iota'$. For simplicity we keep the same kind of notations when we deal with $\E'(\Om)$, $\S'(\R^n)$ and the Colombeau algebras and topological duals connected with them. 

The \emph{generalized singular support} $\singsupp_\g\, T$ of a map $T\in\L(\Gc(\Om),\wt{\C})$ is defined as the complement of the set of points $x$ in $\Om$ such that for some open neighborhood $V$ of $x$ the restriction $T_{\vert_V}$ belongs to $\Ginf(V)$.
\begin{proposition}
\label{prop_comparison_distr}
\leavevmode
\begin{itemize}
\item[(i)] For all $w\in\D'(\Om)$ the maps $\iota_d(w), \iota'(w)\in\L(\Gc(\Om),\wt{\C})$ coincide on $\Gc^\infty(\Om)$;
\item[(ii)] assertion $(i)$ holds with $\D'(\Om)$, $\L(\Gc(\Om),\wt{\C})$ and $\Gc^\infty(\Om)$ replaced by the spaces $\E'(\Om)$, $\L(\G(\Om),\wt{\C})$ and $\Ginf(\Om)$ respectively;
\item[(iii)]assertion $(i)$ holds with $\D'(\Om)$, $\L(\Gc(\Om),\wt{\C})$ and $\Gc^\infty(\Om)$ replaced by the spaces $\S'(\R^n)$, $\L(\GS(\R^n),\wt{\C})$ and $\GSinf(\R^n)$ respectively.
\item[(iv)] For all $w\in\D'(\Om)$
\beq
\label{supp_equality}
\supp\, w\ =\ \supp\, \iota(w)\ =\ \supp\, \iota'(w)\ =\supp\, \iota_d(w)
\eeq
and
\beq
\label{singsupp_equality}
\singsupp\, w\ =\ \singsupp_\g\, \iota(w)\ =\ \singsupp_\g\, \iota'(w)\ = \singsupp_\g\, \iota_d(w).
\eeq
\end{itemize}
\end{proposition}
\begin{proof}
We omit the proof of $(i)$, $(ii)$, $(iii)$ since it can be easily obtained adapting the arguments of Theorem 1.2.63 \cite{GKOS:01} to this context. Concerning $(iv)$ the equality $\supp\, w =\supp\, \iota(w)$ is a feature of the Colombeau algebra $\G(\Om)$ and the definition of support there, while $\supp\, \iota(w)= \supp\, \iota'(w)$ is guaranteed by the fact that the support of a generalized function $u$ does not change when we consider $u$ as an element of $\L(\Gc(\Om),\wt{\C})$ (Proposition 2.13 \cite{GGO:03}). Finally, by definition of $\iota_d$ it is clear that $\supp\, w=\supp\, \iota_d(w)$. Concerning the singular supports we have that $\singsupp\, w =\singsupp_\g\, \iota(w)$ since combining Theorem 25.2 in \cite{O:92} with Proposition 1.2.19 in \cite{GKOS:01} we have that $\iota(w)_{\vert_{\Om'}}\in\Ginf(\Om')$ if and only if $w_{\vert_{\Om'}}\in\Cinf(\Om')$. We now claim that $\singsupp_\g\, u =\singsupp_\g\, (\int_\Om u(y)\cdot\, dy)$ for all $u\in\G(\Om)$. In fact the inclusion $\singsupp_g\, (\int_\Om u(y)\cdot\, dy)\subseteq\singsupp_\g\, u$ is clear and the other one is implied by the embedding $\G(\Om)\subseteq\L(\Gc(\Om),\wt{\C})$. Thus $\singsupp_g\iota(w)=\singsupp_g\iota'(w)$. It remains to prove that $\singsupp\, w=\singsupp_g\iota_d(w)$. Since the inclusion $\singsupp_\g\, \iota_d(w)\subseteq \singsupp\, w$ is clear let us assume that $x_0\in\Om\setminus\singsupp_\g \iota_d(w)$. This implies that there exists a neighborhood $V$ of $x_0$ and a regular generalized function $v\in\Ginf(V)$ such that for all $\chi\in\Cinfc(V)$, for all $\psi\in\Cinfc(V)$ identically $1$ in a neighborhood of the compact support of $\chi$ and for all $u\in\G(\Om)$ we can write
\beq
\label{equality_revised}
\iota_d(w)(\psi u)=\int_\Om v(y)\psi(y)u(y)dy.
\eeq
We recall that $[((\psi w\ast\varphi_\eps)_{\vert_\Om})_\eps]$ and $[((\int_\Om v_\eps(y)\psi(y)\varphi_\eps(x-y)dy)_{\vert_\Om})_\eps]$ belong to $\G(\Om)$ by the embedding of $\D'(\Om)$ into $\G(\Om)$ and by Proposition 2.14 in \cite{GGO:03} respectively. Given $\wt{x}=[(x_\eps)_\eps]\in\wt{\Om}_{\rm{c}}$ we insert $[((\varphi_\eps (x_\eps-\cdot))_{\vert_\Om})_\eps]$ in place of $u$ in \eqref{equality_revised}. From the point value theory it follows that $[((\psi w\ast\varphi_\eps)_{\vert_\Om})_\eps]$ and $[((\int_\Om v_\eps(y)\psi(y)\varphi_\eps(x-y)dy)_{\vert_\Om})_\eps]$ are the same generalized function and then 
\beq
\label{net_revised}
\chi(\psi w\ast\varphi_\eps) -\chi(v_\eps\psi\ast\varphi_\eps) 
\eeq
is a net in $\Neg(\Om)$. Applying $\Delta^m$ to \eqref{net_revised} we have for some constant $c>0$ and for $\eps$ small enough the estimate
\[
\Vert\Delta^m(\chi(\psi w\ast\varphi_\eps))\Vert_{1}\le \Vert\Delta^m(\chi(v_\eps\psi\ast\varphi_\eps))\Vert_1+c.
\] 
which combined with the $\EMinf$-property of $(v_\eps)_\eps$ leads to the following assertion:
\begin{multline}
\label{MO_eq}
\exists N\in\N\, \forall m\in\R\, \exists C>0\, \exists\eta\in(0,1]\, \forall\eps\in(0,\eta]\\  \Vert\Delta^m(\chi(\psi w\ast\varphi_\eps))\Vert_{1}\le C\eps^{-N}.
\end{multline}
By the reasoning employed in the proof of \cite[Theorem 25.2]{O:92}, \eqref{MO_eq} allows to conclude that $\chi w$ is smooth and then $x_0\in\Om\setminus\singsupp\, w$.
\end{proof}
By Remarks \ref{remark_v}, \ref{remark_gtii} we know that even though $\iota_d(w)$ and $\iota'(w)$ coincide on regular generalized functions, they are not the same functional on $\Gc(\Om)$. As an example consider $\delta_0\in\S'(\R^n)\subseteq\Gt(\R^n)$ and note that $\iota'(\delta_0)(u)$ is the class $[(\int_{\R^n}\varphi_\eps(y)u_\eps(y)\, dy)_\eps]$. For a suitable choice of mollifier we have that $\iota_d(\delta_0)\neq \iota'(\delta_0)$ in $\L(\GS(\R^n),\wt{\C})$.

One may reasonably ask if the embedding $\E'(\R^n)\subseteq \S'(\R^n)\subseteq\D'(\R^n)$ can be mimicked by the locally convex $\wt{\C}$-modules $\L(\G(\R^n),\wt{\C})$, $\L(\GS(\R^n),\wt{\C})$ and $\L(\Gc(\R^n),\wt{\C})$, as it is partially done in $(iii)$Theorem \ref{theorem_local_embedding}. We give a negative answer to this question in the following proposition.
\begin{proposition}
\label{prop_embedding_duals}
\leavevmode
\begin{trivlist}
\item[(i)] The map 
\[
\GS(\R^n)\to\G(\R^n):(u_\eps)_\eps+\NS(\R^n)\to (u_\eps)_\eps+\Neg(\R^n)
\]
is not injective.
\item[(ii)] The restriction map 
\[
\L(\GS(\R^n),\wt{\C})\to\L(\Gc(\R^n),\wt{\C}):T\to T_{\vert_{\Gc(\R^n)}}
\]
is not injective. Consequently $\Gc(\R^n)$ is not dense in $\GS(\R^n)$.
\end{trivlist}
\end{proposition}
\begin{proof}
$(i)$ The non-injectivity of the map defined in $(i)$ is due to the fact that we can find a net $(u_\eps)_\eps\in\ES(\R^n)\cap\Neg(\R^n)$ which is not $\S$-negligible. We take a mollifier $\varphi\in\S(\R^n)$ with $\widehat{\varphi}(x_0)(\widehat{\varphi}(x_0)-1)\neq 0$ for some $x_0\in\R^n$. Easy computations show that $(\widehat{\varphi}(\eps x)(\widehat{\varphi}(\eps x)-1))_\eps\in\ES(\R^n)\cap\Neg(\R^n)$ but $(\widehat{\varphi}(\eps x)(\widehat{\varphi}(\eps x)-1))_\eps\not\in\NS(\R^n)$, since by the point value theory induced by $\Gt(\R^n)$ on $\GS(\R^n)$, the net $((\widehat{\varphi}(\eps x)(\widehat{\varphi}(\eps x)-1))_{\vert_{x=x_0/\eps}})_\eps$ is not negligible. Consequently we cannot provide an embedding of $\L(\G(\R^n),\wt{\C})$ into $\L(\GS(\R^n),\wt{\C})$ by mapping $\GS(\R^n)$ into $\G(\R^n)$.
 
$(ii)$ $\Gc(\R^n)$ is a $\wt{\C}$-submodule of $\GS(\R^n)$ and therefore the following restriction map  
\[
\L(\GS(\R^n),\wt{\C})\to\L(\Gc(\R^n),\wt{\C}):T\to T_{\vert_{\Gc(\R^n)}}
\]
is well-defined and continuous. Let us choose a mollifier $\varphi\in\S(\R^n)$ with $\int_{\R^n}\widehat{\varphi}^2(x)(\widehat{\varphi}(x)-1)\, dx\neq 0$. The integral operator $Tu=\int_{\R^n}v(y)u(y)\, dy$ where $v=(\widehat{\varphi}(\eps x)(\widehat{\varphi}(\eps x)-1))_\eps+\NS(\R^n)\in\GS(\R^n)$ belongs to $\L(\GS(\R^n),\wt{\C})$ and since $(\widehat{\varphi}(\eps x)(\widehat{\varphi}(\eps x)-1))_\eps\in\Neg(\R^n)$ its restriction to $\Gc(\R^n)$ is identically $0$. However, $T$ is not $0$ on $\GS(\R^n)$. In fact for $u=(\widehat{\varphi}(\eps x))_\eps+\NS(\R^n)\in\GS(\R^n)$ we have that $Tu=[(\eps^{-n}\int_{\R^n}\widehat{\varphi}^2(x)(\widehat{\varphi}(x)-1)\, dx)_\eps]\neq 0$. In addition this result proves that $\Gc(\R^n)$ is not dense in $\GS(\R^n)$.
\end{proof}
The well-known regularization of distributions via convolution with a mollifier inspires a regularization of $\delta_{\wt{x}}\in\L(\GS(\R^n),\wt{\C})$ via generalized functions in $\GS(\R^n)$.  
\begin{proposition}
\label{prop_regularization}
Let $\rho\in\S(\R^n)$ with $\int_{\R^n}\rho(x)\, dx=1$ and 
\[
v_{\wt{x},q}:=(\rho_{\eps^q}(x_\eps-\cdot))_\eps+\NS(\R^n).
\]
The sequence $(v_{\wt{x},q})_q$ converges to $\delta_{\wt{x}}$ in $\L(\GS(\R^n),\wt{\C})$.
\end{proposition}
\begin{proof}
As in the proof of Proposition \ref{delta_proposition}, $v_{\wt{x},q}$ is a well-defined element of $\GS(\R^n)$. We want to prove that for all bounded subsets $B$ of $\GS(\R^n)$ and for all $p\in\N$ there exists $\overline{q}\in\N$ such that for all $q\ge\overline{q}$ and for all $u\in B$
\beq
\label{est_convergence}
\val_{\wt{\C}}\biggl(\int_{\R^n}v_{\wt{x},q}(y)u(y)\, dy - u(\wt{x})\biggr)\ge p.
\eeq
At the level of representatives we can write
\begin{multline*}
\int_{\R^n}\rho_{\eps^q}(x_\eps-y)u_\eps(y)\, dy - u_\eps(x_\eps)=\int_{\R^n}\rho(y)(u_\eps(x_\eps-\eps^q y)-u_\eps(x_\eps))\, dy\\ =\int_{\R^n}\rho(y)\sum_{|\alpha|=1}\frac{\partial^\alpha u_\eps(x_\eps-\eps^q\theta y)}{\alpha !}(-\eps^q y)^\alpha\, dy
\end{multline*}
and then for $\eps$ small enough $|\int_{\R^n}\rho_{\eps^q}(x_\eps-y)u_\eps(y)\, dy - u_\eps(x_\eps)|\le \eps^{q-N}$ where $N$ depends on the derivatives of order $1$ of $u$. Since $B$ is bounded it is possible to choose the same $N$ for all $u$ in $B$. Hence for every $q\ge p+N$ and $u\in B$ we get \eqref{est_convergence}.
\end{proof}
Note that for any $\wt{x}=[(x_\eps)_\eps]\in\wt{\Om}_{\rm{c}}$ with $x_\eps\in K\Subset\Om$ for $\eps$ small enough, taking a cut-off function $\psi\in\Cinfc(\Om)$ identically 1 in a neighborhood of $K$ the generalized function $v_{\wt{x},q}:=(\psi(\cdot)\rho_{\eps^q}(x_\eps-\cdot))_\eps +\Neg(\Om)$ defines a sequence in $\Gc(\Om)$ which approximates $\delta_{\wt{x}}$ in $\L(\G(\Om),\wt{\C})$.
 
Motivated by Proposition \ref{representation_prop} and the previous regularization of the functional $\delta_{\wt{x}}\in\L(\GS(\R^n),\wt{\C})$ we are able to show the following density result. 
\begin{proposition}
\label{prop_GS_remark}
$\GS(\R^n)$ is dense in $\Gt(\R^n)$.
\end{proposition}
\begin{proof}
The main idea in constructing a sequence of functions in $\GS(\R^n)$ which converges to a given $u\in\Gt(\R^n)$, is to take a regularization of $\delta_x$ for $x$ varying in $\R^n$ and to arrange the representatives in order to have $\S$-moderate estimates. In detail, for $\rho\in\S(\R^n)$ with $\int_{\R^n}\rho(x)dx=1$ we define
\[
u_q:=\biggl(\widehat{\rho}(\eps^q x)\int_{\R^n}\rho_{\eps^q}(x-y)u_\eps(y)\, dy\biggr)_\eps +\NS(\R^n).
\]
$u_q$ makes sense in $\GS(\R^n)$ since for all $\alpha,\beta\in\N^n$
\begin{multline*}
x^\alpha\partial^\beta\biggl(\widehat{\rho}(\eps^q x)\int_{\R^n}\rho_{\eps^q}(x-y)u_\eps(y)\, dy\biggr)\\=\sum_{\gamma\le\beta}\binom{\beta}{\gamma}x^\alpha\partial^\gamma\widehat{\rho}(\eps^q x)\eps^{2q|\gamma|-q|\beta|}\int_{\R^n}\partial^{\beta-\gamma}\rho(y)u_\eps(x-\eps^q y)\, dy
\end{multline*}
and therefore $$\displaystyle\sup_{x\in\R^n}\biggl|x^\alpha\partial^\beta\big(\widehat{\rho}(\eps^q x)\int_{\R^n}\rho_{\eps^q}(x-y)u_\eps(y)\, dy\big)\biggr|=O(\eps^{-q|\alpha+\beta|-qN-M})$$ if $\sup_{x\in\R^n}(1+|x|)^{-N}|u_\eps(x)|=O(\eps^{-M})$. In particular $u_q$ is uniquely determined by $u$ once we have chosen the mollifier.

Since the topology on $\Gt(\R^n)$ is based on $\GtS(\R^n)=\cap_{m\in\N}\cup_{N\in\N}\G^m_{N,\S}(\R^n)$ we preliminarily show that if $m\ge 1$ and  $u\in\G^m_{N,\S}(\R^n)$ then $u_q\to u$ in $\G^{m-1}_{N+1,\S}(\R^n)$. It is convenient to write $\widehat{\rho}(\eps^q x)\int_{\R^n}\rho_{\eps^q}(x-y)u_\eps(y)\, dy - u_\eps(x)$ as the sum of $s_{1,\eps}(x):=\widehat{\rho}(\eps^q x)\int_{\R^n}\rho(y)(u_\eps(x-\eps^q y)-u_\eps(x))\, dy$ and $s_{2,\eps}(x):=(\widehat{\rho}(\eps^q x)-1)u_\eps(x)$. For every derivative $\alpha$ with $|\alpha|\le m-1$ and $(u_\eps)_\eps\in\E^m_{N}(\R^n)$ with $\sup_{x\in\R^n,|\gamma|\le m}(1+|x|)^{-N}|\partial^\gamma u_\eps(x)|=O(\eps^{-M})$, $|\partial^\alpha s_{1,\eps}(x)|$ is estimated by 
\begin{multline}
\label{s_1}
\sum_{\beta\le\alpha}c_{\alpha,\beta}\eps^{q|\beta|}|\partial^\beta\widehat{\rho}(\eps^q x)|\int_{\R^n}|\rho(y)|\sum_{|\gamma|=1}\frac{|\partial^{\alpha-\beta+\gamma}u_\eps(x-\eps^q\theta y)|}{\gamma !}|(-\eps^q y)^\gamma|\, dy\\ \le c'\eps^{q-M}(1+|x|)^N
\end{multline}
while
\begin{multline}
\label{s_2}
|\partial^\alpha s_{2,\eps}(x)|\le\sum_{\beta\le\alpha}c_{\alpha,\beta}\sum_{|\gamma|=1}\biggl|\frac{\partial^{\beta+\gamma}\widehat{\rho}(\eps^q\theta x)}{\gamma !}(\eps^qx)^\gamma\biggr||\partial^{\alpha-\beta}u_\eps(x)|\\
\le c\eps^{q-M}(1+|x|)^{N+1}.
\end{multline}
\eqref{s_1} and \eqref{s_2} together yield the convergence of $u_q$ to $u$ in $\G_{N+1,\S}^{m-1}(\R^n)$.

At this point, noting that for all $m\ge 1$ each $u\in\GtS(\R^n)$ is an element of $\G^{m}_{N,\S}(\R^n)\subseteq \G^{m-1}_{\tau,\S}(\R^n)$ for some $N\in\N$, combining the previous arguments with the inductive limit topology on $\G^{m-1}_{\tau,\S}(\R^n)$, we have that $u_q$ tends to $u$ in $\G^{m-1}_{\tau,\S}(\R^n)$. It follows  that the sequence $(u_q)_q$ converges to $u$ in $\GtS(\R^n)$. Finally we consider the map $\jmath:\GtS(\R^n)\to \Gt(\R^n):(u_\eps)_\eps+\NS(\R^n)\to (u_\eps)_\eps+\Nt(\R^n)$. For any $u\in\Gt(\R^n)$, representatives $(u_\eps)_\eps$ of $u$ and $(u_\eps)_\eps+\NS(\R^n)\in\jmath^{-1}(u)$, the sequence $(u_q)_q\subseteq \GS(\R^n)$ constructed as above has limit $(u_\eps)_\eps+\NS(\R^n)$ in $\GtS(\R^n)$. This implies $\jmath(u_q)\to \jmath((u_\eps)_\eps+\NS(\R^n))=u$ in $\Gt(\R^n)$ as desired.
\end{proof}
\begin{remark}
\label{remark_final_chapter4}
From the proof of the previous proposition it is also clear that $\GS(\R^n)$ is dense in $\GtS(\R^n)$. Hence by restriction we obtain the embeddings of $\L(\GtS(\R^n),\wt{\C})$ and $\L(\Gt(\R^n),\wt{\C})$ into $\L(\GS(\R^n),\wt{\C})$ respectively. These embeddings are continuous if we equip the duals with the topology $\beta_b$ of uniform convergence on bounded subsets.
\end{remark}

\bibliographystyle{abbrv}
%\bibliographystyle{abbrv}
%\bibliography{claudia}

\newcommand{\SortNoop}[1]{}

\end{document}